%% file: art023.tex
\documentclass[11,reqno,oneside]{amsart}
\usepackage[margin=3cm]{geometry}

\usepackage{subfigure}
\usepackage{tikz}
\usetikzlibrary{calc}
\usepackage[boxed,oldcommands]{algorithm2e}
\usepackage{hyperref}
\usepackage{tabularx}
\usepackage{afterpage}
\usepackage{lscape}
\usepackage{multirow}
\usepackage{multicol}
\usepackage{arydshln}
\usepackage{cancel}
\usepackage{makecell}
\usepackage{amsfonts,amscd}
\usepackage{amssymb,amsmath}
\usepackage{tikz}

\usepackage[utf8]{inputenc}
\usepackage[T1]{fontenc}
\usepackage{lmodern}
\usepackage{nicefrac}

\usepackage{autonum}
\usepackage[numbers,sort&compress]{natbib}
\usepackage{stmaryrd,xparse}
\usepackage{doi}
\usepackage{epstopdf}

\include{definitions}

\def\tfc{\boldsymbol{\Phi}}
\def\cfunc{\boldsymbol{e}}

\def\Eq#1{(\ref{#1})}
\def\x{{\boldsymbol{x}}}

\def\VbddcC{\tilde{\mathbb{V}}_{\rm c}}

\def\VbddcF{{\tilde{\mathbb{V}}_{\rm f}}}
\def\VbddcF1{{\tilde{\mathbb{V}}_{\rm f}}}
\def\VbddcC1{{\tilde{\mathbb{V}}_{\rm c}}}

\def\u{\boldsymbol{u}}
\def\y{\boldsymbol{y}}

\def\shortTitle{Nonlinear parallel-in-time multilevel Schur complement solvers for ordinary differential equations}

\def\myKeywords{Time parallelism, ordinary differential equations, domain decomposition, nonlinear solver, scalability}

\begin{document}

\title[\shortTitle]{Nonlinear parallel-in-time Schur complement solvers for ordinary differential equations}

\author[S. Badia]{ Santiago Badia$^\dag$}

\thanks{$\dag$ Universitat Polit\`ecnica de Catalunya, Jordi Girona1-3, Edifici C1, E-08034 Barcelona $\&$ Centre Internacional de M\`etodes Num\`erics en Enginyeria, Parc Mediterrani de la Tecnologia, Esteve Terrades 5, E-08860 Castelldefels, Spain E-mail: {\tt sbadia@cimne.upc.edu}. SB gratefully acknowledges the support received from the Catalan Government through the ICREA Acad\`emia Research Program. This work has been partially funded by the project MTM2014-60713-P from the ``Ministerio de Econom\'ia, industria y Competitivad'' of Spain.}

\author[Marc Olm]{Marc Olm$^\ddag$}
\thanks{$\ddag$ Universitat Polit\`ecnica de Catalunya, Jordi Girona1-3, Edifici C1, E-08034 Barcelona $\&$ Centre Internacional de M\`etodes Num\`erics en Enginyeria, Parc Mediterrani de la Tecnologia, Esteve Terrades 5, E-08860 Castelldefels, Spain E-mail: {\tt molm@cimne.upc.edu}. MO gratefully acknowledges the support from the Ag\`encia de Gesti\'o i d'Ajuts Universitaris i de Recerca, under the FI-AGAUR 2015 grant. }
 
\date{\today}

\begin{abstract}
  In this work, we propose a parallel-in-time solver for linear and nonlinear ordinary differential equations. The approach is based on an efficient multilevel solver of the Schur complement related to a multilevel time partition. For linear problems, the scheme leads to a fast direct method.  Next, two different strategies for solving nonlinear ODEs are proposed. First, we consider a Newton method over the global nonlinear ODE, using the multilevel Schur complement solver at every nonlinear iteration. Second, we state the global nonlinear problem in terms of the nonlinear Schur complement (at an arbitrary level), and perform nonlinear iterations over it.  Numerical experiments show that the proposed schemes are weakly scalable, i.e., we can efficiently exploit increasing computational resources to solve for more time steps the same problem. 
\end{abstract}

\maketitle



\tableofcontents



\section{Introduction} 

{
\setalcaphskip{1ex}
\setlength{\algomargin}{1em}
\restylealgo{ruled} 
\linesnumbered
\dontprintsemicolon
\Setnlsty{textrm}{ }{:}
\SetFuncSty{\texbf}
\SetKwComment{comment}{}{}{}
\SetCommentSty{textrm}


At the beginning of the next decade supercomputers are expected to reach a peak performance of one exaflop/s, which implies a 100 times improvement with respect to current supercomputers. This improvement will not be based on faster processors, but on a much larger number of processors (in a broad sense). This situation will certainly have an impact in large scale computational science and engineering (CSE). Parallel algorithms will be required to exhibit much higher levels of concurrency, keeping good scalability properties.



When dealing with transient problems, since information always moves forward in time, one can exploit sequentiality. However, the tremendous amounts of parallelism to be exploited in the near future certainly motivates to change this paradigm. One of the motivations to exploit higher levels of parallelism will be to reduce the time-to-solution. In the simulation of ordinary differential equations (ODEs), the way to go is to exploit concurrency in time. The idea is to develop parallel-in-time solvers that provide the solution at all time values in one shot, instead of the traditional sequential approach that exploits the arrow of time. If scalable parallel-in-time solvers are available, the use of higher levels of parallelism will certainly reduce the time-to-solution.

Parallel-in-time solvers are receiving rapidly increasing attention. Different iterative methods have been considered so far, e.g., the  \emph{parareal} method \cite{parareal} or spectral deferred-correction time integrators \cite{pfasst}. With regard to direct methods, time-parallel methods can be found in \cite{Gander_50years}. In general these methods can exploit low levels of concurrency \cite{christlieb_parallel_2010} or are tailored for particular types of equations \cite{gander_paraexp:_2013}. We refer to \cite{Gander_50years} for an excellent up-to-date review of time parallelism. It has also motivated the development of space-time parallel solvers for transient partial differential equations (PDEs) \cite{falgout_parallel_2014,badia_space-time_2017}. 

In this work, we propose a parallel-in-time solver for ODEs that relies on the well-known Schur complement method in linear algebra. For linear (systems of) ODEs, the approach can be understood as a Schur complement solver in time. 
 When the coarse problem is too large compared to the local problems, due to the structure of the coarse problem, we can consider recursively the Schur complement strategy, leading to multilevel implementations, in order to push forward scalability limits. The method can be applied to $\theta$-methods, discontinuous Galerkin (DG) methods, Runge-Kutta methods, and BDF methods. We also note that the proposed method can also be understood as a parareal scheme in which the coarse solver is automatically computed in such a way that the scheme is a direct method (convergence in one iteration is assured). (The interpretation of the parareal method as an approximation of the Schur complement problem has already been pointed out in \cite{falgout_parallel_2014}.) As a result, the proposed method solves the drawback of the parareal scheme, i.e., its poor parallel efficiency, inversely proportional to the number of iterations being required by the iterative algorithm. One of the messages of this work is to show that the approximation of the Schur complement in parareal methods does not really pay the price when (just by roughly multiplying by two the number of operations) one can have a highly scalable direct parallel-in-time solver.

In order to extend these ideas to nonlinear PDEs, we consider two different strategies. First, we consider a global linearization of the problem in time using, e.g., Newton's method. We note that the idea of a global linearization of nonlinear ODEs to exploit time-parallelism is not new. It was already considered by Bellen and Zenaro in 1989 \cite{bellen_parallel_1989}. (In any case, the solvers proposed in  \cite{bellen_parallel_1989} are different from the ones presented herein. They are restricted to the Steffensen's linearization method, which leads to a diagonal problem per nonlinear iterations, where parallelization can obviously be used.) After the linearization of the problem, we consider the Schur complement solver commented above at every nonlinear iteration. A second strategy consists in applying the nonlinear Schur complement strategy first, and next to consider the linearization of such operator, leading to nested nonlinear iterations.

The parallel-in-time ideas in this work can naturally be blended with domain decomposition ideas to design highly scalable space-time parallel solvers. We have combined these ideas with a multilevel balancing DD by constraints (BDDC) preconditioner (see \cite{mandel_multispace_2008,tu_three-level_2007}) in \cite{badia_space-time_2017}, by judiciously choosing the quantities to be continuous among processors in a space-time partition, i.e. \cite{dohrmann_preconditioner_2003}. The resulting space-time parallel solver has been proved to be scalable on thousands of processors for different transient (non)linear PDEs.

The outline of the article is as follows. In Sect. \ref{sec:prob-st}, we state the problem. We introduce a time-parallel direct method for linear ODEs based on the computation of a multilevel time Schur complement in Sect. \ref{sec:ODEs}. In Sect. \ref{sec:nonlinear}, we extend the method to nonlinear ODEs, by combining first a Newton linearization step with a Schur complement linear solver and next considering nonlinear Schur complement problems. We present a detailed set of numerical experiments in Sect. \ref{sec:numerical_experiments}, showing the excellent scalability properties of the proposed methods. Finally, we draw some conclusions in Sect. \ref{sec:conclusions}.

\section{Statement of the problem}\label{sec:prob-st}
In this section, we develop a parallel direct solver for the numerical approximation of ordinary differential equations (ODEs). We consider a system of ODEs of size $m_{\rm unk}$: 
\def\u{\boldsymbol{u}}
\def\x{\boldsymbol{v}}
\def\f{\boldsymbol{f}}
\def\e{\boldsymbol{e}}
\def\v{\boldsymbol{v}}
\def\y{\boldsymbol{y}}
  \def\f{\boldsymbol{f}}
  \def\g{\boldsymbol{g}}
  \def\kappaop{\boldsymbol{\kappa}}
\begin{align}\label{eq-ode}
  \displaystyle\frac{{\rm d}\u(t)}{{\rm d}t} + \kappaop (t,\u(t)) = \boldsymbol{0}, \qquad  \u(t^0) = \u_0, 
\end{align}
for $t \in  (t^0 = 0,T]$. Let us assume that $\kappaop(\cdot,\cdot)$ is continuous with respect to the first argument and Lipschitz continuous with respect to the second argument, and that existence and uniqueness holds.

\def\ie{k}
\def\elt{{\delta_\ie}}
\def\elsbt{{\delta_\is^\ie}}
\def\sbt{{\Delta_\is}}
\def\ns{N}
\def\ne{K}
\def\is{n}
For the time interval $[0,T]$, we define a hierarhical multilevel partition as follows (see Fig. \ref{fig-part} for a detailed illustration). We define a (level-0) time partition $\{0 = t^0_0 , t^1_0, \ldots, t^{n_0}_0 = T\}$ into $n_0$ time elements. Next, we consider a (level-1) coarser time partition $\{0 = t^0_1, \ldots, t_1^{n_1} = T \}$ into $n_1$ time subdomains (or level-1 elements), defined by aggregation of elements at the previous level, i.e., for every $i \in \{ 0,\ldots,n_1 \}$ there exists an $m(i) \in \{ 0, \ldots, n_0\}$ such that $t_1^i = t_0^{m(i)}$. We proceed recursively, creating coarser partitions for higher levels. We define the time element $i$ at level-$k$ as the time interval $(t_k^i,t_k^{i+1})$.
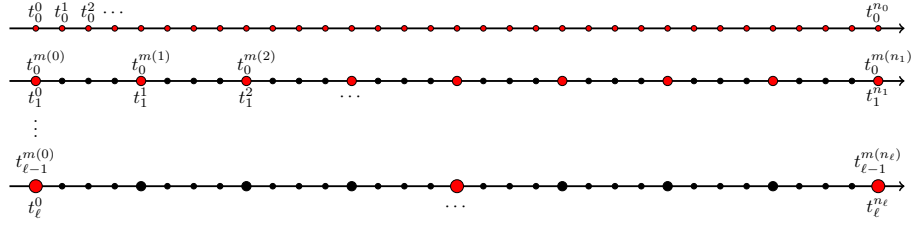
\begin{figure}
  \begin{center}
\scalebox{0.7}{\input{fig_time_part}}
\end{center}
\caption{Multilevel time partition of $[0,T]$. The subindex in $t_\alpha^\beta$ denotes the partition level, whereas the superindex denotes the time step in such partition. In order to relate time values of two constitutive levels, we use the notation $t_\alpha^\beta = t_{\alpha-1}^{m(\beta)}$, i.e., the time value $t_\alpha^\beta$ corresponds to the time step $m(\beta)$ at the previous level.}\label{fig-part}
\end{figure}

\def\Ac{\boldsymbol{\mathcal{A}}}
\def\Ec{\boldsymbol{\mathcal{E}}}
\def\Kc{\boldsymbol{\mathcal{K}}}
\def\Jc{\boldsymbol{\mathcal{J}}}
We will present the method in a general way that is independent of the time integration scheme being used. We define the nonlinear operators $\Ac_0^{i+1} : \u^i \mapsto \u^{i+1}$ for $i \in \{ 0, \ldots, n_0-1\}$, such that, given the initial value { $\u^{i}$}, solves \Eq{eq-ode} in $(t^i_0,t^{i+1}_{0})$, and provides $\u^{i+1} = \u(t^{i+1}_0)$. We conceptually state the solver at the continuous level, even though one can consider different time integration schemes instead, e.g., $\theta$-methods, DG methods, or Runge-Kutta methods. The use of BDF-type schemes requires some further ellaboration. We refer to Sect. \ref{subsec:other_time_int} for more details.

\section{An ODE direct solver}\label{sec:ODEs}

\def\Km{\mathbf{K}}
\def\Em{\mathbf{E}}
\def\Fm{\mathbf{F}}
\def\uv{\mathbf{u}}
\def\ev{\mathbf{e}}
\def\fv{\mathbf{f}}
\def\gv{\mathbf{g}}
\def\yv{\mathbf{y}}
\def\vv{\mathbf{v}}
\def\zv{\mathbf{z}}

\def\tfc{\boldsymbol{\Phi}}
  \def\Ic{\boldsymbol{\mathcal{I}}}
  \def\Fc{\tfc}
  \def\Kc{\boldsymbol{\mathcal{K}}}
  \def\U{\boldsymbol{U}}
  \def\G{\boldsymbol{G}}
  In this section, we assume that $\kappaop(t,\cdot)$ is a linear operator. (The nonlinear extension is described in Sect. \ref{sec:nonlinear}.) In this case, it is easy to check that $\Ac_0^{i+1}(\cdot)$ is an affine mapping, and we have $\Ac_0^{i+1}(\u^i) = \Fc_0^{i+1} \u^i + \g_0^{i+1}$, where $\Fc_0^{i+1}$ is a linear operator (an $m_{\rm unk} \times m_{\rm unk}$ matrix). Both $\Fc_0^{i+1}$ and $\g_0^{i+1}$ can be explicitly computed from $\kappaop(\cdot,\cdot)$ and $\f$ for a given time integration scheme. Thus, the global problem \Eq{eq-ode} for linear ODEs can be stated in algebraic form as:
\begin{align} \label{global_ODE}
\left( \begin{array}{cccc}
\Ic          &         &         &        \\
-\Fc_0^{1}         & \Ic    &         &        \\
           & \ddots  &  \ddots &        \\
&         &  -\Fc_0^{n_0}     &  \Ic  
\end{array} \right)
\left( \begin{array}{c}
\u^0_0   \\
\u^1_0    \\
 \vdots \\
\u^{n_0}_0
\end{array} \right)
=
\left( \begin{array}{c}
\u_0   \\
\g_0^1    \\
 \vdots \\
\g_0^{n_0}
\end{array} \right),
\end{align} where $\Ic$ is the identity matrix (of size
$m_{\rm unk} \times m_{\rm unk}$). We can also represent the global system \Eq{global_ODE} in compact notation with
\begin{align}
  \Km_0 \uv_0 = \gv_0, \quad \hbox{or equivalently} \quad 
\left( \begin{array}{cc}
\Km_0^{II} & \Km_0^{I\Gamma}  \\
\Km_0^{\Gamma I} & \Km_0^{\Gamma \Gamma}      
\end{array} \right)
\left( \begin{array}{c}
\uv^I_0   \\
\uv^\Gamma_0  
\end{array} \right)
=
\left( \begin{array}{c}
\gv_0^I    \\
\gv_0^{\Gamma}
       \end{array}\right),
\end{align}
where we have considered a segregation of degrees of freedom (DOFs) at level-0 $\uv_0$ into the \emph{interface} DOFs $\uv_0^\Gamma$, namely the time step values that are also in the level-1 partition, and the \emph{interior} dofs $\uv_0^I$. $\Km_0$ is a 2-banded lower block-triangular matrix. In order to define the Schur complement problem, we first consider an \emph{interior correction} of the problem at hand,
  \begin{align}\label{int_corr}
  \Km_0^{II}  \vv_0^I = \gv_0^I,
  \end{align}
  i.e., we solve the system  \Eq{global_ODE} in the subspace of vectors that vanish in the level-1 time values $t_1^0, \ldots, t_1^{n_1}$. The solution of \Eq{int_corr} involves $n_1$ independent local ODE problems: solve for $i=0,\ldots,n_1-1$
\begin{align}\label{local_ODE}
\left( \begin{array}{cccc}
\Ic          &         &         &        \\
-\Fc_0^{m(i)+1}         & \Ic    &         &        \\
           & \ddots  &  \ddots &        \\
&         &  -\Fc_0^{m(i+1)-1}     &  \Ic 
\end{array} \right)
\left( \begin{array}{c}
\x^{m(i)}_{0}   \\
\x^{m(i)+1}_{0}    \\
 \vdots \\
\x^{m(i+1)-1}_{0}
\end{array} \right)
=
\left( \begin{array}{c}
0  \\
\g_0^{m(i)+1}    \\
 \vdots \\
\g_0^{m(i+1)-1}
\end{array} \right).
\end{align}
After the interior correction, we must solve the problem
\begin{align}\label{interface_correction}
\left( \begin{array}{cc}
\Km_0^{II} & \Km_0^{I\Gamma}  \\
\Km_0^{\Gamma I} & \Km_0^{\Gamma \Gamma}      
\end{array} \right)
\left( \begin{array}{c}
\delta \uv_0^{I}   \\
\uv^\Gamma_0  
\end{array} \right)
=
\left( \begin{array}{c}
         \mathbf{0}
         \\
\gv_0^{\Gamma}
       \end{array}\right), \quad \hbox{and compute} \quad \uv_0^I = \vv_0^I + \delta \uv_0^I.
  \end{align}
In order to solve \Eq{interface_correction}, we define the following extension operator (usually denoted as the harmonic extension operator in the frame of domain decomposition solvers for PDEs):
  \begin{align}\label{extension}
\left( \begin{array}{cc}
\Km_0^{II} & \Km_0^{I\Gamma}  \\
\mathbf{0} & \mathbf{I}_\Gamma  
\end{array} \right)
\left( \begin{array}{c}
\Em_0^I   \\
\Em_0^\Gamma  
\end{array} \right)
=
\left( \begin{array}{c}
\mathbf{0}    \\
\mathbf{I}_\Gamma  
       \end{array}\right), \quad \hbox{thus} \quad
    \Em_0 = 
\left( \begin{array}{c}
- (\Km_0^{II})^{-1} \Km_0^{I\Gamma}    \\
\mathbf{I}_\Gamma  
       \end{array}\right).
  \end{align}
  Thus, using the fact that $ \uv^\Gamma_0 = \uv_1$, the Schur complement of level-0 reads:
  \begin{align}\label{schur_1} 
(\Km_0^{\Gamma \Gamma} +  \Km_0^{I \Gamma} \Em_0^{I}) \uv_1 = \gv_0^{\Gamma} - \Km_0^{I\Gamma} \vv_0^{I}, \quad \hbox{represented by} \quad \Km_1 \uv_1 = \gv_1.
  \end{align}
 The Schur complement of the level-0 system is the level-1 problem. By construction, the extension operator solution of \Eq{extension} can be written as a block-diagonal matrix  $\Em_0 = {\rm diag}(\ev_{(0)}, \ev_{(1)}, \cdots, \ev_{(n_1-1)},\Ic)$, which involves $n_1 \times {m_{\rm unk}}$ independent local ODE problems: solve for $i=0,\ldots,n_1-1$
\begin{align}\label{coarse_func}
\left( \begin{array}{cccc}
\Ic          &         &         &        \\
-\Fc_0^{m(i)+1}         & \Ic    &         &        \\
           & \ddots  &  \ddots &        \\
&         &  -\Fc_0^{m(i+1)-1}     &  \Ic 
\end{array} \right)
\left( \begin{array}{c}
(\cfunc_{(i)})^{m(i)}_{0}   \\
(\cfunc_{(i)})^{m(i)+1}_{0}    \\
 \vdots \\
(\cfunc_{(i)})^{m(i+1)-1}_{0}
\end{array} \right)
=
\left( \begin{array}{c}
\Ic  \\
0    \\
 \vdots \\
0
\end{array} \right).
\end{align}
After some manipulation, the Schur complement problem { \Eq{schur_1}} can be written as:
\begin{align}\label{global_Schur}
\left( \begin{array}{cccc}
\Ic          &         &         &        \\
-\tfc^{1}_1        & \Ic    &         &        \\
           & \ddots  &  \ddots &        \\
&         &  -\tfc^{n_1}_1     &  \Ic 
\end{array} \right)
\left( \begin{array}{c}
\u^0_1   \\
\u^1_1    \\
 \vdots \\
\u^{n_1}_1
\end{array} \right)
=
\left( \begin{array}{c}
\u_0   \\
\g_1^{m(1)}    \\
 \vdots \\
\g_1^{m(n_1)}
       \end{array} \right),
  \end{align}
  where $\tfc^{i}_1 \doteq \tfc_0^{m(i)} \cfunc_{(i-1)}^{m(i)-1}$ and $\g_1^i \doteq \g_0^{m(i)} + \tfc_0^{m(i)-1} \v_0^{m(i)-1}$. Thus, the Schur complement matrix $\Km_1$ at the next level has also the same structure as the original problem, i.e., it is a ODE-type solver over the coarser level-1 time partition. The computation of the interior correction, the extension operator $\Em_0$, the Schur complement matrix $\Km_1$, and the Schur complement right-hand side $\gv_1$ can readily be computed in parallel. In Alg. \ref{alg:ass_schur} (for $k = 0$), all the steps to assemble the Schur complement problem are listed, indicating on the left-hand side the line task corresponding level.
{
\setalcaphskip{1ex}
\setlength{\algomargin}{1em}
\restylealgo{ruled} 
\linesnumbered
\dontprintsemicolon
\Setnlsty{textrm}{ }{:}
\SetFuncSty{\texbf}
\SetKwComment{comment}{}{}{}
\SetCommentSty{textrm}

\begin{algorithm}
  \small{
  \KwData{$\Km_{k}$, $\gv_{k}$}
  \KwResult{$\vv_{k}$, $\Em_{k}$, $\Km_{k+1}$,$\gv_{k+1}$}
  Compute the interior correction $\vv_{k}$ solution of \Eq{int_corr}, by solving the $n_{k+1}$ local ODE problems \Eq{local_ODE} for $i=0, \ldots,n_{k+1}-1$  \comment*[r]{{$k-1$}} 
  Compute the extension operator $\Em_{k}$ solution of \Eq{extension}, by solving the $n_{k+1} \times m_{\rm unk}$ problems \Eq{coarse_func} for $i=0, \ldots,n_{k+1}-1$, and compute $\Km_k^{I \Gamma} \Em_0^I$ and $-\Km_k^{I\Gamma} \uv_k^I$ in \Eq{schur_1} \comment*[r]{{$k-1$}} 
    Assemble the Schur complement system \Eq{global_Schur}, i.e., $\Km_{k+1}$ and $\gv_{k+1}$  (see \Eq{schur_1}) \comment*[r]{${k-1} \rightarrow k$} 
\caption{Schur complement set-up (level-$k$) \label{alg:ass_schur}}}
\end{algorithm}
}
In a two-level implementation, i.e., $\ell = 1$, the Schur complement problem \Eq{interface_correction} would be computed in serial. It would finally lead to $\uv_{0} = \vv_{0} + \Em_{0} \uv_1$. This approach leads to $n_1$ indepedent level-0 problems of size $\frac{n_0}{n_1}$ and one coarse level-1 problem of size $n_1$. Clearly, when $n_1$ increases, the coarse problem becomes the bottleneck of the simulations.  In order to push forward the scalability limits of this approach, we can consider a multilevel Schur complement technique. When $n_1$ exceeds $\frac{n_0}{n_1}$, since \Eq{global_Schur} has the same structure as the original system \Eq{global_ODE}, one can consider the same Schur complement approach for the level-1 system, leading to a three-level algorithm. We can proceed recursively to include an arbitrary number of levels. In Alg. \ref{alg:ml_schur}, we state the multilevel Schur complement ODE solver. We comment on the parallel efficiency and computational cost of this algorithm in Sect. \ref{sec:par_ef}. 
{
\setalcaphskip{1ex}
\setlength{\algomargin}{1em}
\restylealgo{ruled} 
\linesnumbered
\dontprintsemicolon
\Setnlsty{textrm}{ }{:}
\SetFuncSty{\texbf}
\SetKwComment{comment}{}{}{}
\SetCommentSty{textrm}

\begin{algorithm}
  \small{
  \KwData{$\Km_0$, $\gv_0$}
  \KwResult{$\uv_0 = \Km_0^{-1} \gv_0$}
  \For{$k = 0, \ldots,\ell-1$}
  {
    Call Alg. \ref{alg:ass_schur} with $(\Km_k,\gv_k)$ to get $(\vv_k,\Em_{k},\Km_{k+1},\gv_{k+1})$ \comment*[r]{$k,k+1$}
  }
  Solve $\Km_\ell \uv_\ell = \gv_\ell$ \comment*[r]{$\ell$}
  \For{$k = \ell-1, \ldots,0$}
  {  
    Compute $\uv_k = \vv_k + \Em_{k} \uv_{k+1}$ \comment*[r]{$k$}
  }
\caption{Multilevel Schur complement ODE solver  \label{alg:ml_schur}}}
\end{algorithm}
}

\subsection{Application to different time integrators}\label{subsec:other_time_int}

The previous approach can straightforwardly be used for $\theta$-methods. For DG and Runge-Kutta methods, we can require multiple intermediate stages to move from $t_0^{i}$ to $t_0^{i+1}$. In the DG case, we have some additional time values per time element. We can consider that all the element values but the last one are eliminated at the element level, using the so-called static condensation technique. In this case, the resulting discrete problem can be stated as in \ref{global_ODE}. The DG method being used only affects the expression of $\Fc_0^{i+1}$ and $\g_0^{i+1}$. (We note that DG methods do not satisfy time causality at the element level, but it does not affect the lower 2-banded block-triangular structure after the elimination of ``interior'' element values.) We proceed analogously for multi-stage Runge-Kutta methods.

BDF schemes (of second and higher order) slightly differ from the fact that the computation of the solution at a given time step not only requires value from the previous time step, but some additional stages. For a BDF(X) scheme, the system matrix in \Eq{global_ODE} is a (X+1)-banded lower block-triangular matrix. It affects the concept of interface nodes $\Gamma$; to decouple the global problems into local problems, we require to increase the size of the interface X times. As a result, the coarse-scale space dimension is X times larger, as well as the number of coarse space basis functions being computed. In this sense, high order BDF schemes are be a bad choice when dealing with parallel computations, since the interface among subdomains increases, with the corresponding computational cost, due to a loss of locality with respect to the continuous problem. High order Runge-Kutta or DG methods are better suited for time-parallel computations. In any case, after considering the modification described above, the techniques proposed in this work can be applied to BDF methods.

\subsection{Parareal interpretation}

The multilevel Schur complement solver  defined in Alg.~\ref{alg:ml_schur} can also be understood as a (multilevel) parareal scheme. In the parareal scheme, we consider a coarse solver to provide initial conditions to local fine solvers. Instead, in the Schur complement method, one first computes the (fine) interior correction, which is required to the assembly of the right-hand side in the coarse solver. The method above can be stated in a different way, by defining the restriction operator $\Fm_0$ as follows:
\begin{align}
  \left(
  \begin{array}{cc}
     \Fm_0^I   &    \Fm_0^\Gamma
  \end{array} \right)            
  \left( \begin{array}{cc}
\Km_0^{II} & 0 \\
\Km_0^{\Gamma I} & \mathbf{I}_\Gamma  
\end{array} \right)
=
\left( \begin{array}{c}
0    \\
\mathbf{I}_\Gamma  
       \end{array}\right), \quad \hbox{thus} \quad
    \Fm_0 = 
\left( \begin{array}{cc}
-  \Km_0^{\Gamma I} (\Km_0^{II})^{-1}  & \mathbf{I}_\Gamma
       \end{array}\right).
\end{align}
Analogously to the extension operator $\Em_0$, the computation of the restriction operator is a block-diagonal matrix  $\Fm_0 = {\rm diag}(\Ic,\fv_{(0)}^T, \fv_{(1)}^T, \cdots, \fv_{(n_1-1)}^T)$, which involves $n_1 \times m_{\rm unk}$ independent local (backwards) ODE problems:  solve for $i=0,\ldots,n_1-1$
\def\trfc{\boldsymbol{\Psi}}
\begin{align}
\left( \begin{array}{cccc}
\Ic          &   -(\Fc_0^{m(i)+2})^T       &         &        \\
        & \Ic    &   \ddots      &        \\ 
           &   &  \ddots &       -(\Fc_0^{m(i+1)})^T   \\
&         &      &  \Ic 
\end{array} \right)
\left( \begin{array}{c}
\f^{m(i)+1}_{(i)}   \\
\f^{m(i)+2}_{(i)}    \\
 \vdots \\
\f^{m(i+1)}_{(i)}
\end{array} \right)
=
\left( \begin{array}{c}
0  \\
0    \\
 \vdots \\
\Ic
\end{array} \right).
\end{align}
The well-posedness of the backward problem is a direct consequence of the well-posedness of its transpose, the forward problem. Thus, the Schur complement problem reads: 
$$
\Km_1 \uv_1 = \gv_0^\Gamma + \Fm_0^I \gv_0^I = \gv_1.
$$
The coarse parareal problem in the Schur complement method is of Petrov-Galerkin type, using as coarse trial space the range of $\Em_0$ and as coarse test space the range of $\Fm_0^T$, i.e,
\begin{align}\label{alt_schur}
\Km_1 = \Fm_0 \Km_0 \Em_0, \quad \gv_1 = \Fm_0 \gv_0.
\end{align}
It is easy to check that { \Eq{schur_1}} and  \Eq{alt_schur} are equivalent. The fine solver is simply the interior correction \Eq{int_corr}. (We note that fine and coarse corrections are independent, since they are $\Km_0$-orthogonal.). The definition of the coarse space is automatic and the method is not an iterative but a direct solver. As a result, the method does not suffer from the low parallel efficiency of parareal methods, which is proportional to the inverse of parareal iterations \cite{parareal}. Even though the implementation that involves the computation of the trial and test coarse spaces is the one being used in non-symmetric PDE solvers with inexact Schur complement preconditioning (see, e.g., \cite{badia_implementation_2013}), it is not convenient for ODE solvers, since it involves an additional fine solver.

Fig. \ref{fig-toy_dG_coarse_phis} depicts forward and backward coarse functions for a test problem. The local oscillations in the DG(1) and DG(2) schemes are due to the fact that time causality does not hold inside the element. (We note that the previous developments have been considered after eliminating (using element-wise solvers) all the element values but the last one (in time).)  Fig. \ref{fig-toy_dG} shows coarse level-1, fine level-0, and full solution for a selected simple problem with different coarse DOFs position consideration.

\begin{figure}[h!!]
\begin{center}
\begin{tabular}{c@{}c@{}c@{}}
DG(0) & DG(1) & DG(2) \\
\includegraphics[width=0.32\textwidth]{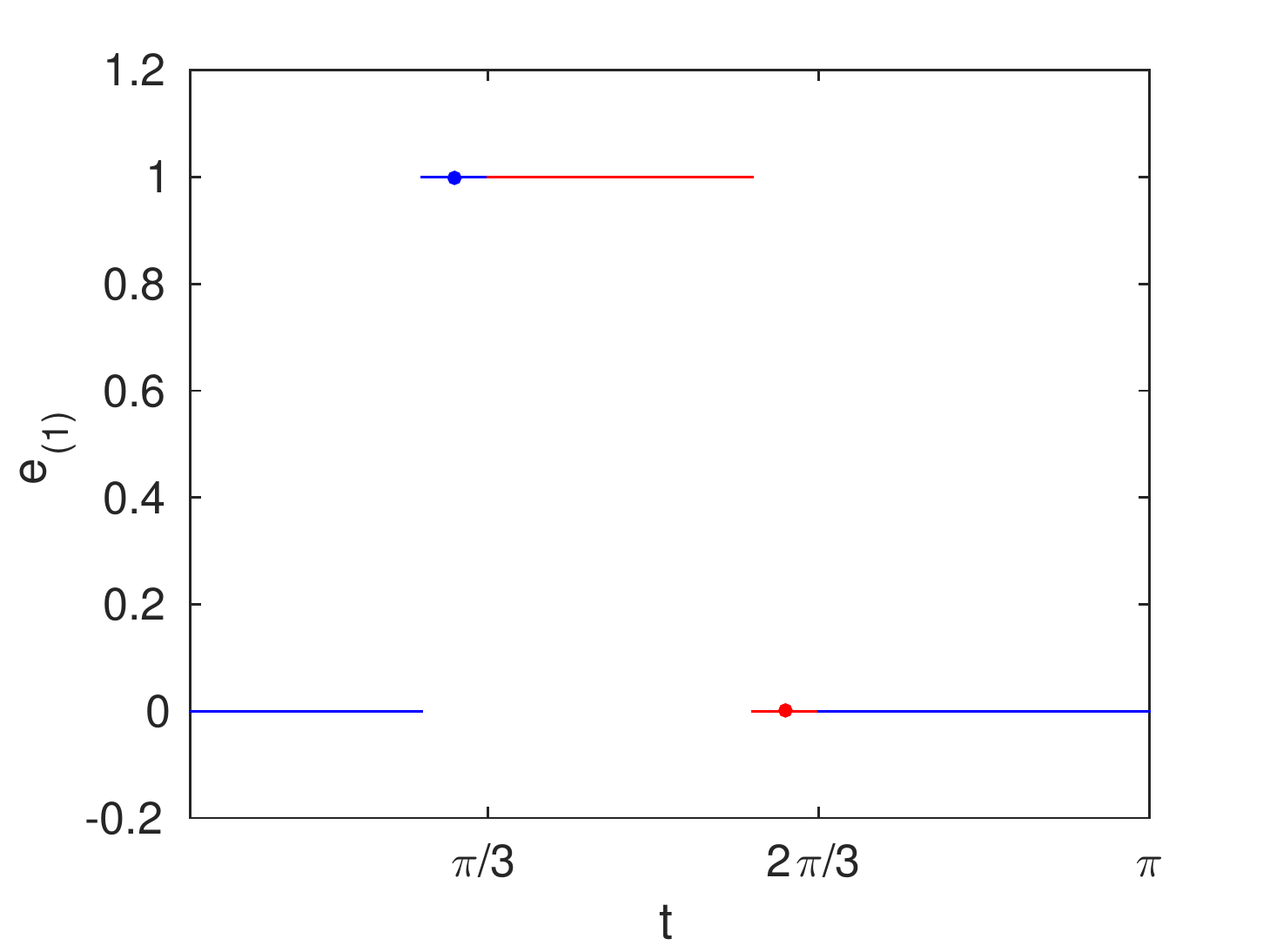}          &
\includegraphics[width=0.32\textwidth]{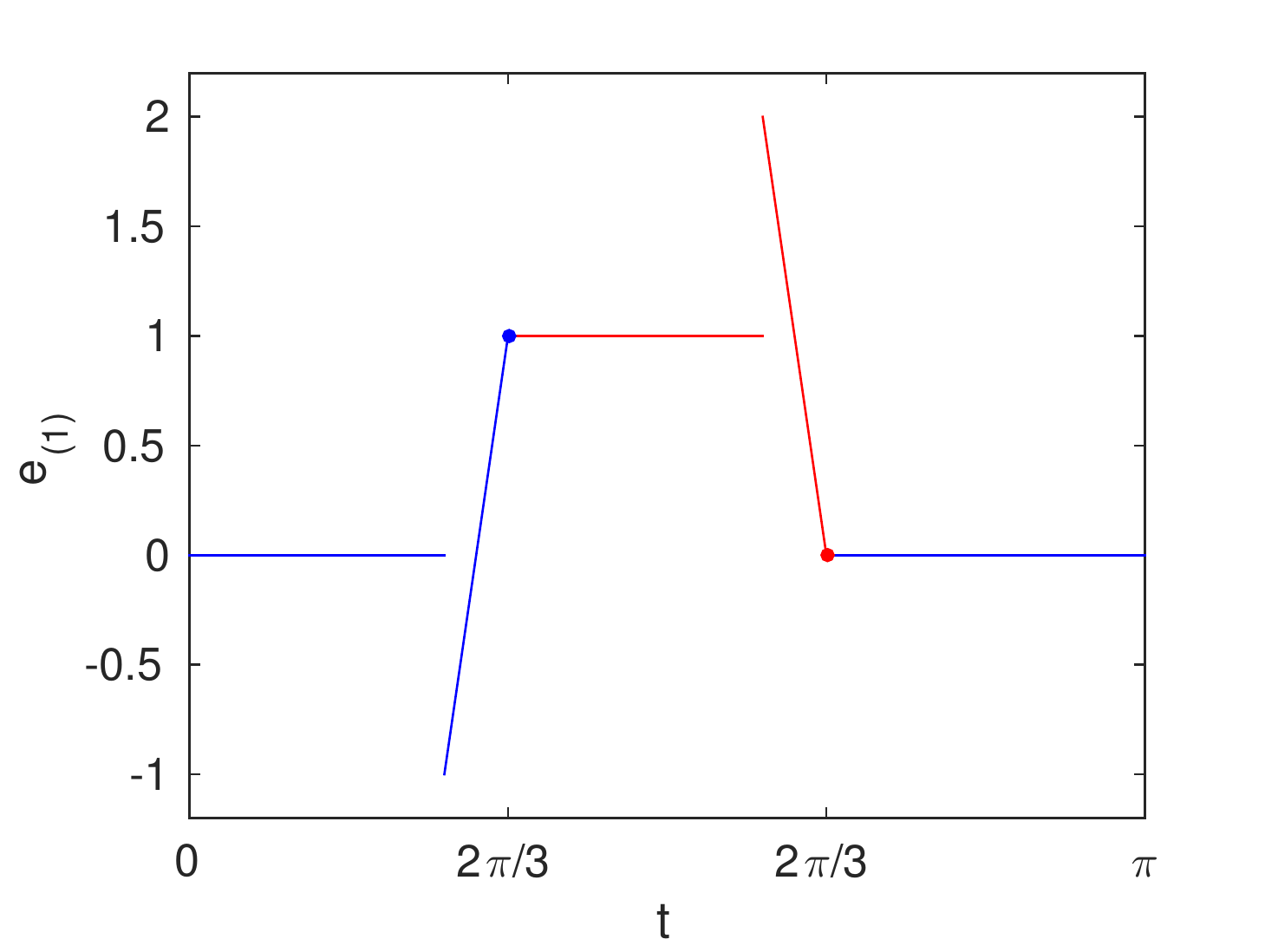}      &
\includegraphics[width=0.32\textwidth]{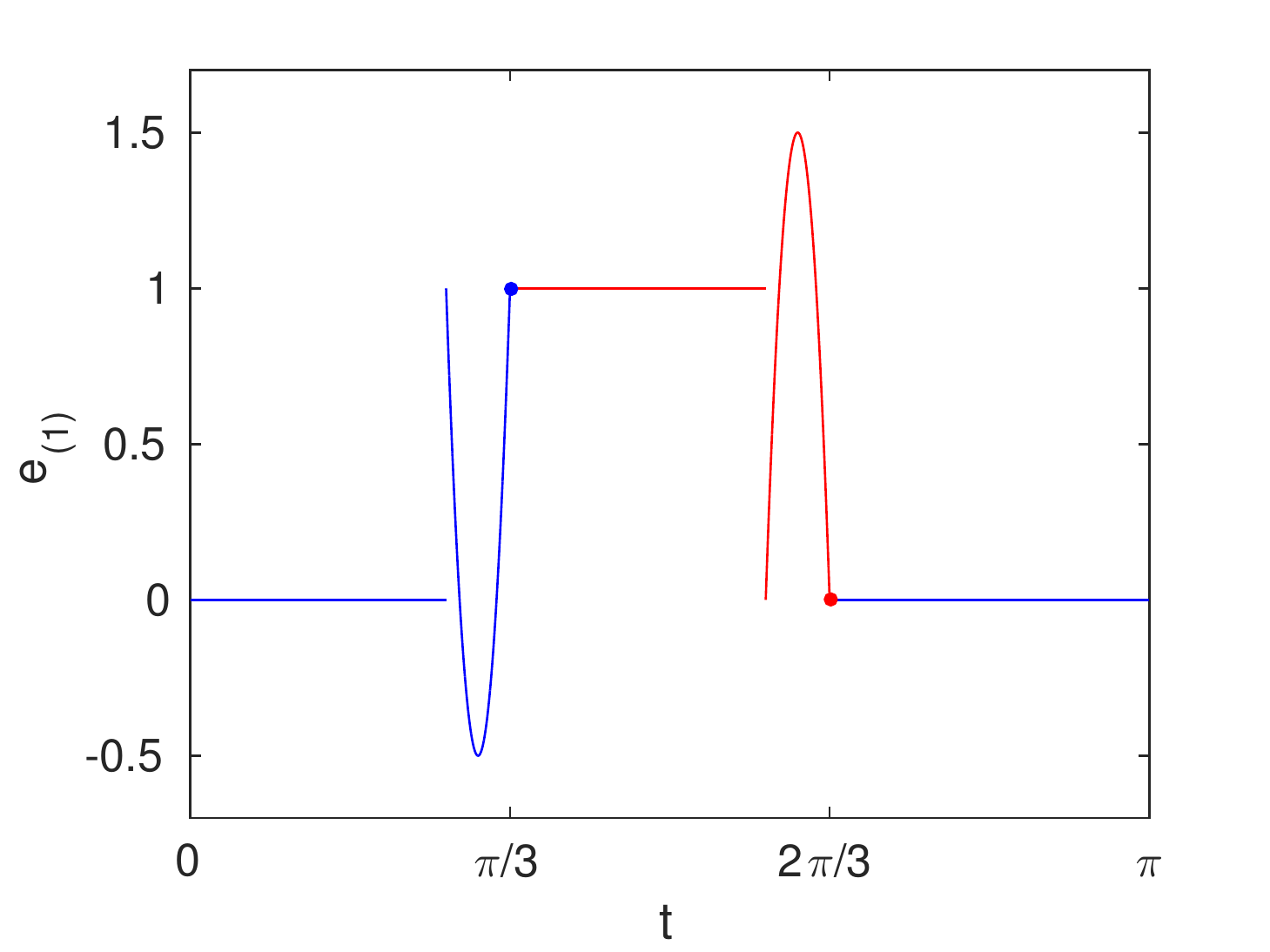}      \\
\multicolumn{3}{c}{(a) \hspace{0.5cm} Coarse shape function $\ev_{(1)}$ for $t_1^1 = \pi/3$}   \\
\includegraphics[width=0.32\textwidth]{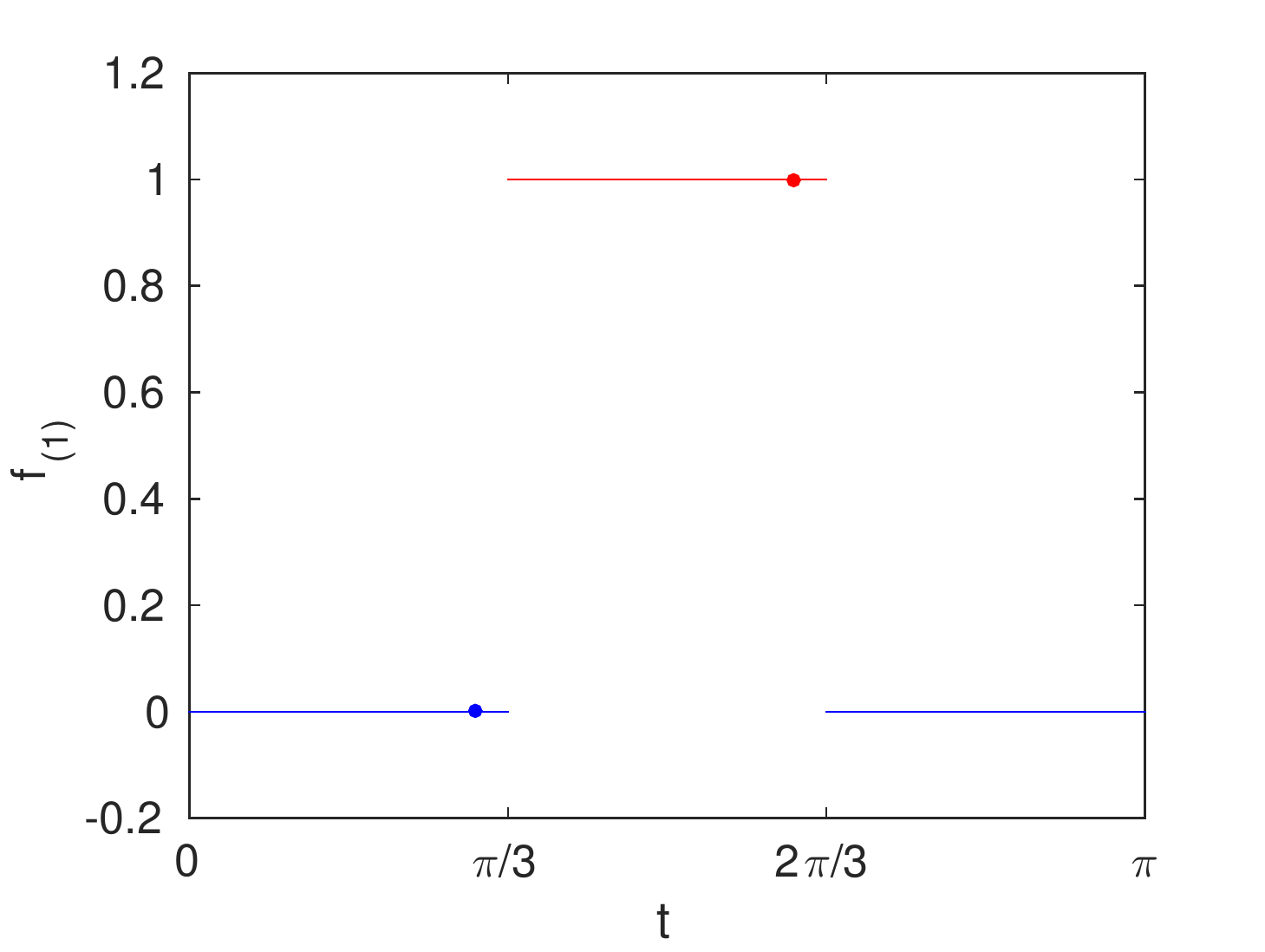}          &
\includegraphics[width=0.32\textwidth]{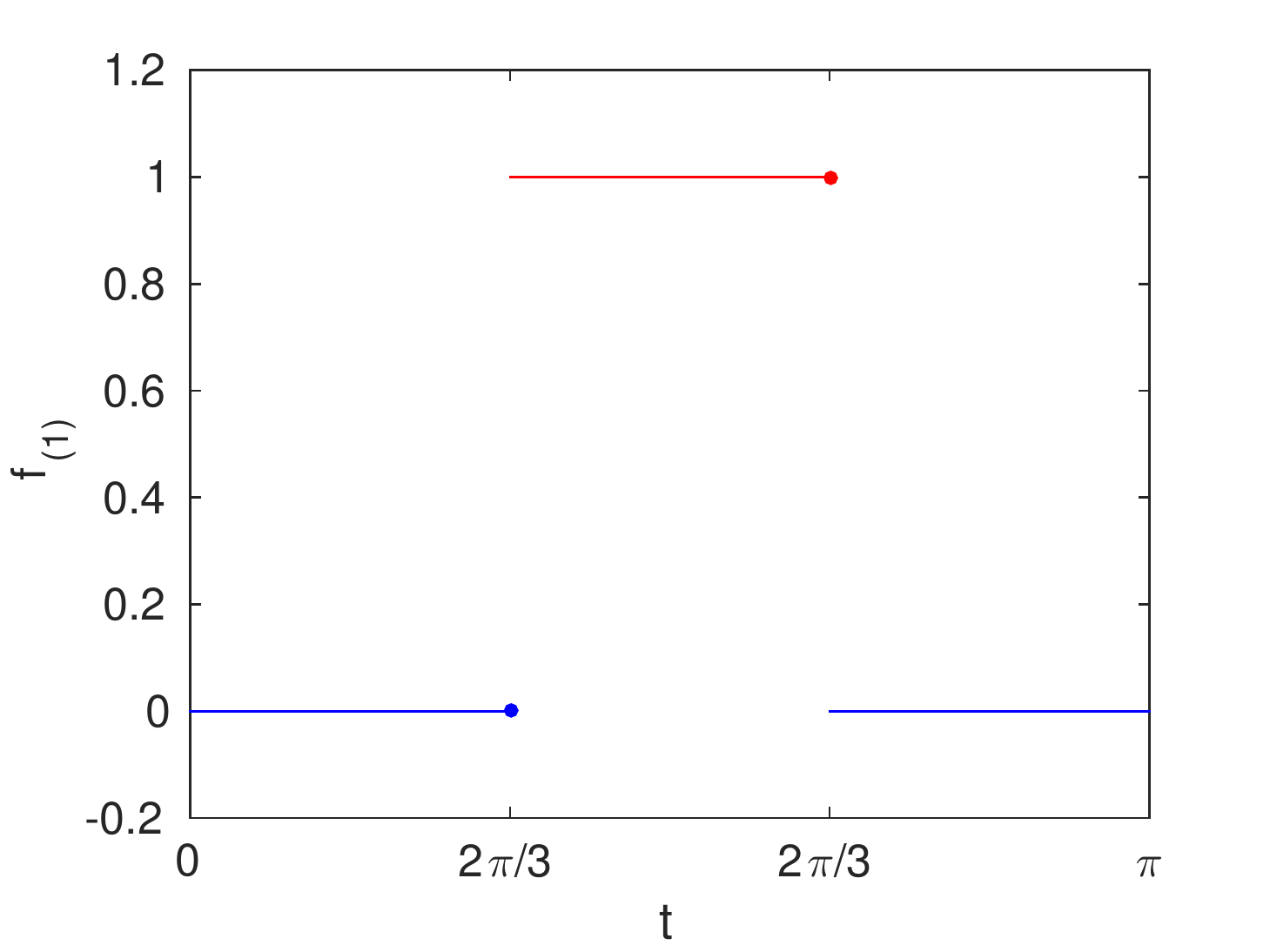}     &
\includegraphics[width=0.32\textwidth]{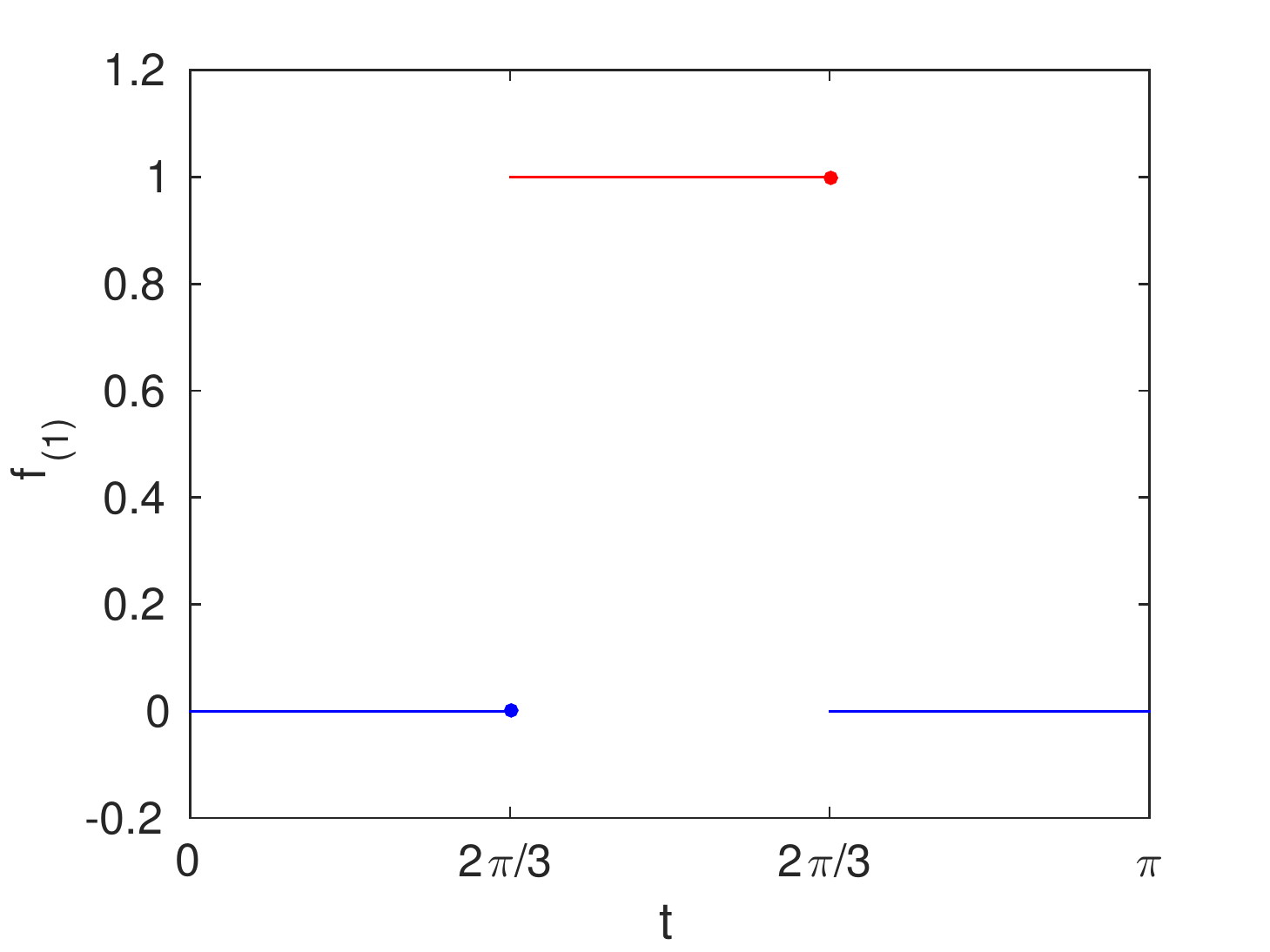}     \\
\multicolumn{3}{c}{(b) \hspace{0.5cm} Coarse shape function $\fv_{(1)}$ for $t_1^1 = \pi/3$}  \\
\end{tabular}
\end{center}
\caption{Coarse shape functions $\ev_{(1)}$ and  $\fv_{(1)}$ ($t_1^1 = \pi/3$) for the simple transport operator $\frac{du}{dt}$.  A partition of an interval $[0,\pi]$ into three subdomains is considered. Each subdomain is partitioned into 5 time elements. Sub-intervals are depicted in consecutive different colour in order to aid visualization, and dots aid to identify coarse DOFs position (in the center of the element for DG(0)).}
\label{fig-toy_dG_coarse_phis}
\end{figure} 

\begin{figure}[h!!]
\begin{center}
\begin{tabular}{c@{}c@{}c@{}}
DG(0) & DG(1) & DG(2) \\
\includegraphics[width=0.32\textwidth]{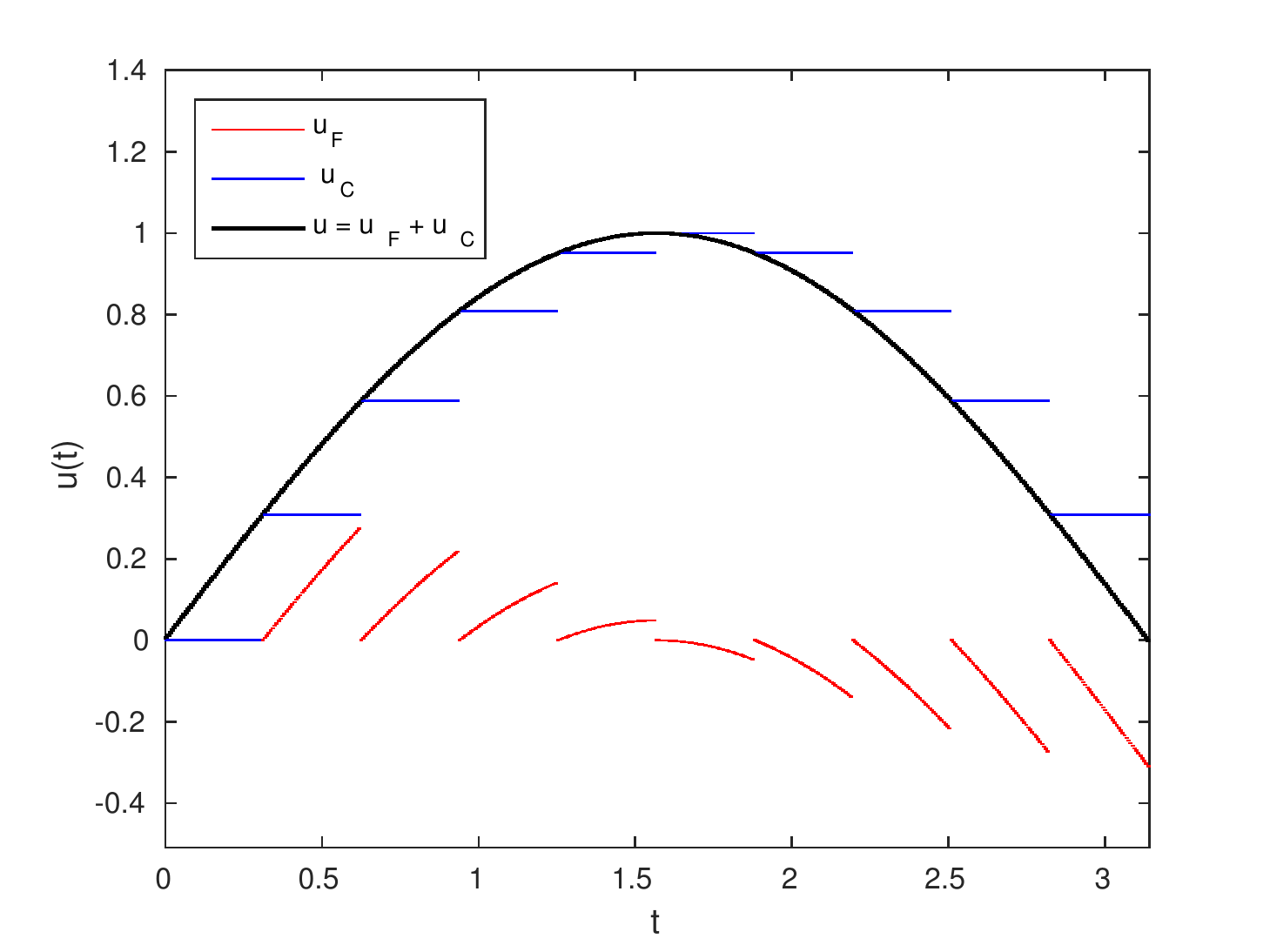}                    &
\includegraphics[width=0.32\textwidth]{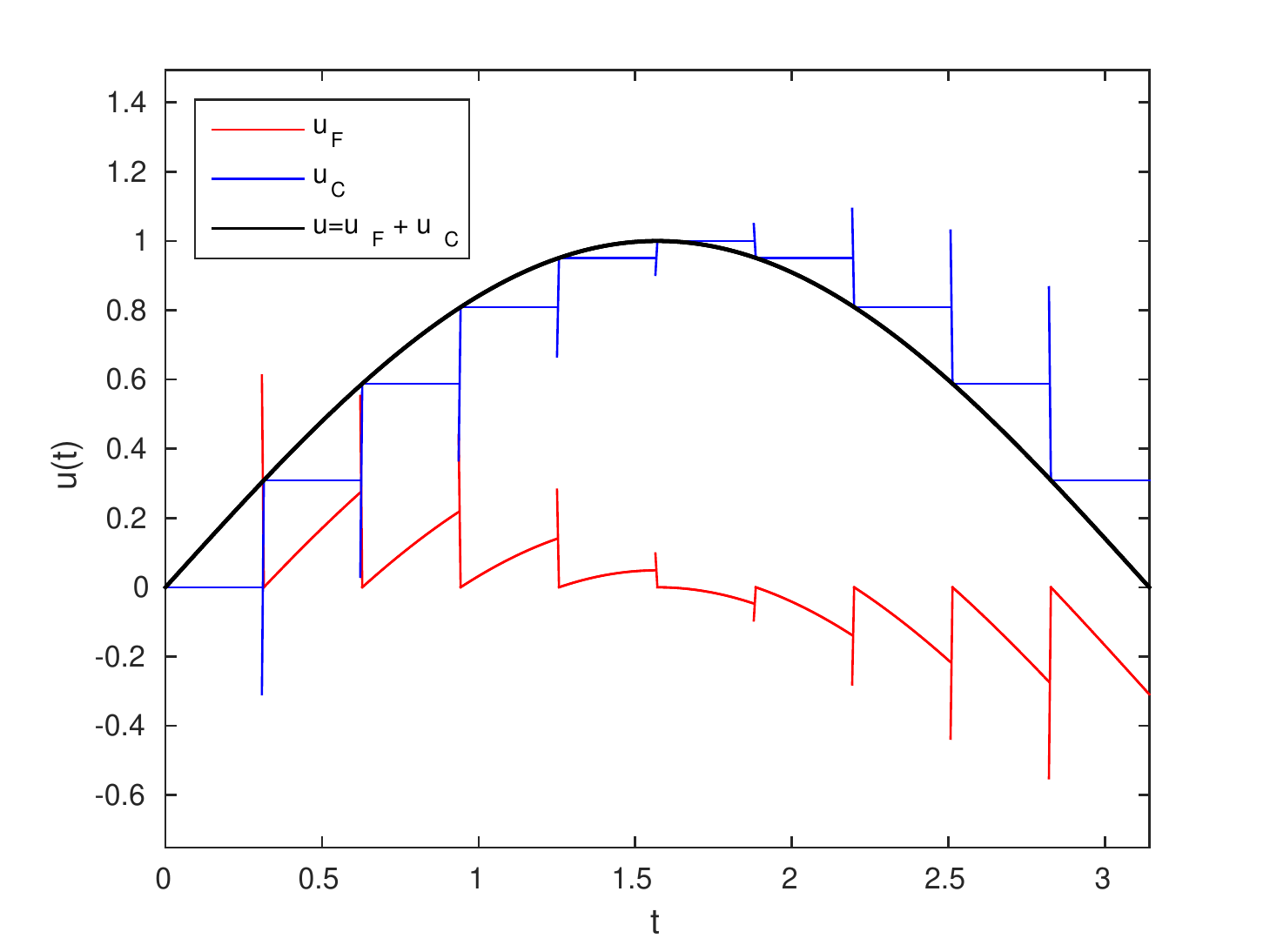}                   &
\includegraphics[width=0.32\textwidth]{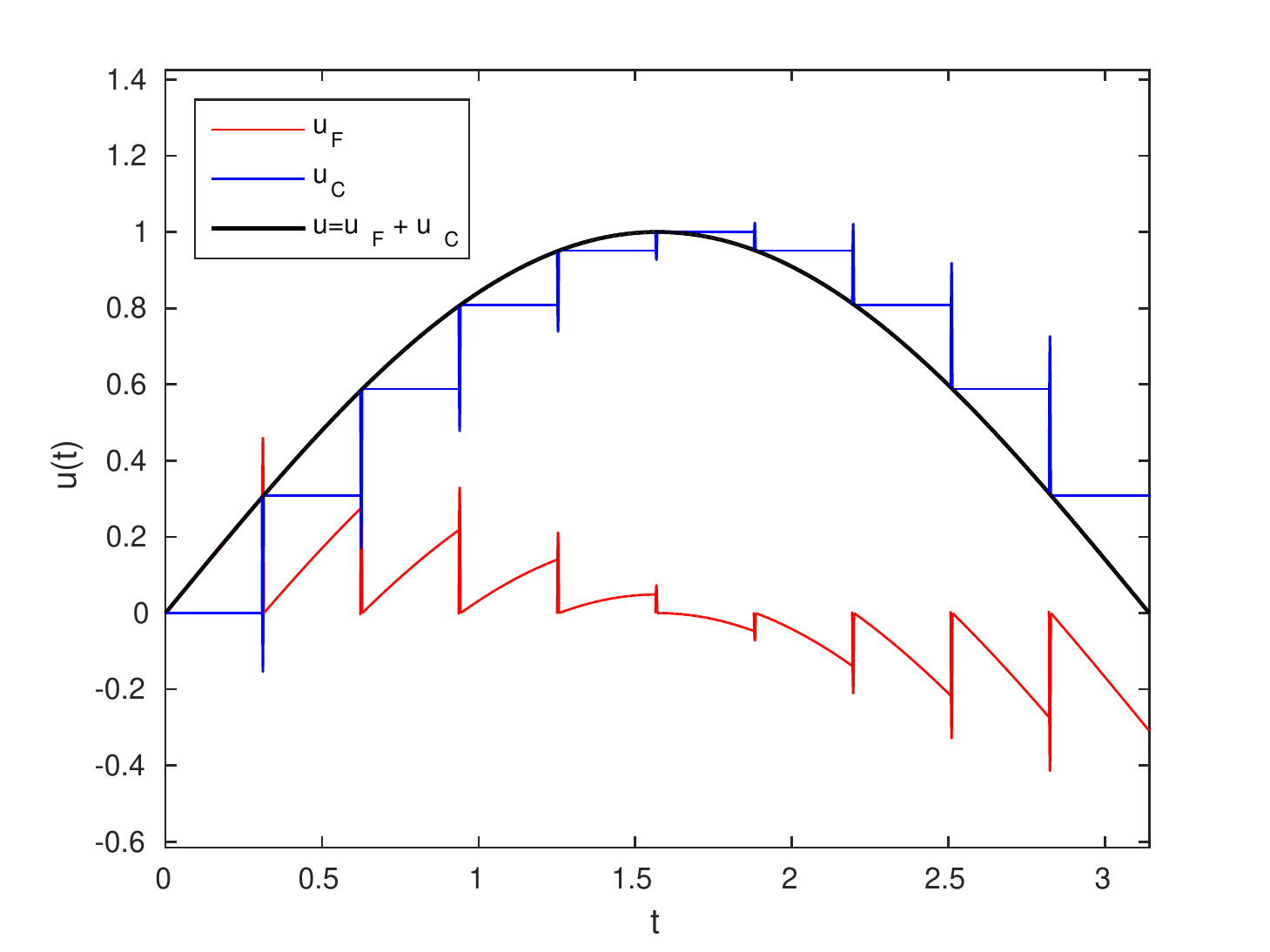}          
\end{tabular}
\end{center}
\caption{ Decomposition of $\uv_0$ into fine $\Em_0 \uv_0$ and coarse $\uv_1$ component for the simple problem $\partial_t u = \cos(t)$, on $[0,\pi]$. The time domain is partitioned into 10 subdomains and each subdomain is an aggregation of 50 elements.  }
\label{fig-toy_dG}\end{figure}

\subsection{Parallel efficiency}\label{sec:par_ef}

Let us consider the multilevel Schur complement solver in Alg. \ref{alg:ml_schur} for a linear ODE. Using the same approach as in multigrid methods, we can consider a fixed coarsening ratio $\theta$ between subdomains, and leave free the number of levels $\ell$ required for a particular simulation with $n_0$ time step values. The number of processors being used is assumed to be equal to $n_1 = n_0/\theta$, i.e., the number of subdomains at level-1, and the number of time steps per processor at all levels is $\theta$ (at the last level it can be smaller). Thus, we define $\ell + 1 = {\rm ceiling} (\frac{\log{n_0}}{\log{\theta}})$, where ${\rm ceiling}(A)$ returns the least integer greater than or equal to $A$.

Alg. \ref{alg:ml_schur} requires at  levels $0, \ldots, \ell-1$ to solve $1+m_{\rm unk}$ local linear ODE problems with $\theta$ time steps per processor (see Alg. \ref{alg:ass_schur}), one to compute the interior correction and $m_{\rm unk}$ to compute the extension operator. One linear ODE with at most $\theta$ time steps must be solved at the last level. The different levels must be computed in a sequential way.

Based on the { solver} described above, we can estimate the order of floating-point operations (FLOPs) required to solve the linear ODE. Sequentially, it is of the order of ${\rm FLOP_0} \approx n_0(m_{\rm unk}^2 + m_{\rm unk})$, since we require $m_{\rm unk}^2 + m_{\rm unk}$ operations per time step. On the other hand, the number of FLOPs required to compute the same problem in parallel using Alg. \ref{alg:ml_schur} is the one needed to solve $(1+m_{\rm unk})$ problems of size $n_0$, $n_0/\theta$, $\ldots$, $n_0/\theta^\ell$. { Using the geometric sequence sum formula, we get ${\rm FLOP}_p \approx  {\rm FLOP_0} (1+m_{\rm unk}) (1-\frac{1}{\theta^{\ell+1}})(1 - \frac{1}{\theta})^{-1} < {\rm FLOP_0} (1+m_{\rm unk}) (1 - \frac{1}{\theta})^{-1}$,} but these operations can be performed in parallel exploiting distributed memory machines. With regard to time, the total CPU time of the multilevel Schur complement in Alg. \ref{alg:ml_schur}  is the aggregation of the CPU time of the solution of $(1+m_{\rm unk})$ linear ODEs with $\theta$ time steps for all levels. Thus, the parallel CPU time is ${\rm CPU}_p \approx  \ell \theta (m_{\rm unk}^2 + m_{\rm unk})(1+m_{\rm unk})$, whereas the serial CPU time is ${\rm CPU}_0 \approx n_0 (m_{\rm unk}^2 + m_{\rm unk})$. As a result, the CPU time linearly depends on the size of the global system in a very mild logarithmic way, i.e., it increases with $\log{n_0}$ due to the expression of $\ell$, and the parallel algorithm will rapidly lead to a shorter time-to-solution than the serial solver.\footnote{We note that the communications are not considered in the previous estimations. In any case, it has been experimentally observed that the communication time in (incomplete) multilevel Schur complement methods is small compared to the computation time up to almost half a million tasks in \cite{badia_multilevel_2016}.} As a result, the speed-up of the proposed direct solver is $S = {\rm CPU}_0/{\rm CPU}_{p} \approx {\rm P} \ell^{-1} (1+m_{\rm unk})^{-1}$. The speed-up is quasi-linear (it only increases with $\log{n_0}$), and thus algorithmically strongly scalable. Analogously, the method is algorithmically weakly scalable, because the total CPU time does not depend on the number of processors or global system size (appart from the logarithmic term).

\section{Iterative solvers for nonlinear ODEs}\label{sec:nonlinear}

In this section, we assume that $\kappaop(\cdot,\cdot)$ is nonlinear. The nonlinear ODE \Eq{eq-ode} in $(0,T]$ can be stated in compact form as $\Ac_0(\uv_0) = \mathbf{0}$, which can also be split into interior and interface time steps as follows:
\begin{align}\label{nonlinear_split}
\left(
\begin{array}{c}
  \Ac_0^I ( \uv_0  ) \\
  \Ac_0^\Gamma ( \uv_0 )
\end{array}
\right) = \mathbf{0}.
\end{align} We consider two different types of solvers for the nonlinear problem.

\subsection{Newton-Schur complement methods}
In order to solve the nonlinear ODE, we can use Newton's method over the global-in-time problem, and solve at every iteration a linear ODE using the multilevel Schur complement in Alg. \ref{alg:ml_schur}. 
We can compute the Jacobian matrix related to \Eq{nonlinear_split} around a point $\bar\uv_0$ as:
\begin{align}\label{linearization_ode}
\Jc_0 (\bar\uv_0) \doteq \frac{\partial \Ac_0  }{\partial \uv_0}(\bar\uv_0) = \frac{ {\rm d} \bar \uv_0}{ {\rm d} t} + \frac{\partial \Kc_0}{\partial \uv_0}(\bar\uv_0),
\end{align}
where we have used the fact that $\Ac_0$ has two terms, one related to the time derivative and the other one related to $\kappaop(t,\cdot)$ (see \ref{eq-ode}). Thus, to solve \Eq{linearization_ode} involves the solution of a linear ODE of the form \Eq{global_ODE}.

The resulting multilevel Newton-Schur complement solver for nonlinear ODEs is stated in Alg. \ref{alg:newton_schur}. Figure \ref{fig-nonlinear_ode} shows the solution iterates using this algorithm  for a selected nonlinear  problem with known analytic solution $u = \rm{sin}(t)$. Each solution update obtained from a linearized problem is solved with the direct solver in Alg. \ref{alg:ml_schur} using two levels, i.e., $\ell = 1$. The nonlinear iterations of the NEwton-Schur complement methods do not depend on the partition being used, since it is a linearization of the global problem and a (parallel) direct solver is used at every nonlinear iteration. In any case, the comparison against the sequential approach is more complicated here, since different nonlinear iterations (local vs. global) are being used in every case. We note that we can readily consider other iterative methods, e.g., Picard's method or Anderson acceleration techniques. In the case of Picard's method, the computation of the Jacobian is not required. Instead, the operator $\Ac_0$ must be written as $\Ac_0(\uv_0) = {\widetilde{\Ac}_0(\uv_0)}\uv_0$ and use at every nonlinear iteration the linear operator ${\widetilde{\Ac}_0(\bar\uv_0)}$ instead of $\Jc_0 (\bar\uv_0)$. In Sect. \ref{sec:numerical_experiments}, we consider a hybrid Picard-Newton nonlinear solver. 
\begin{figure}[h!!]
\centering
\subfigure[Convergence history]{\includegraphics[width=0.3\textwidth]{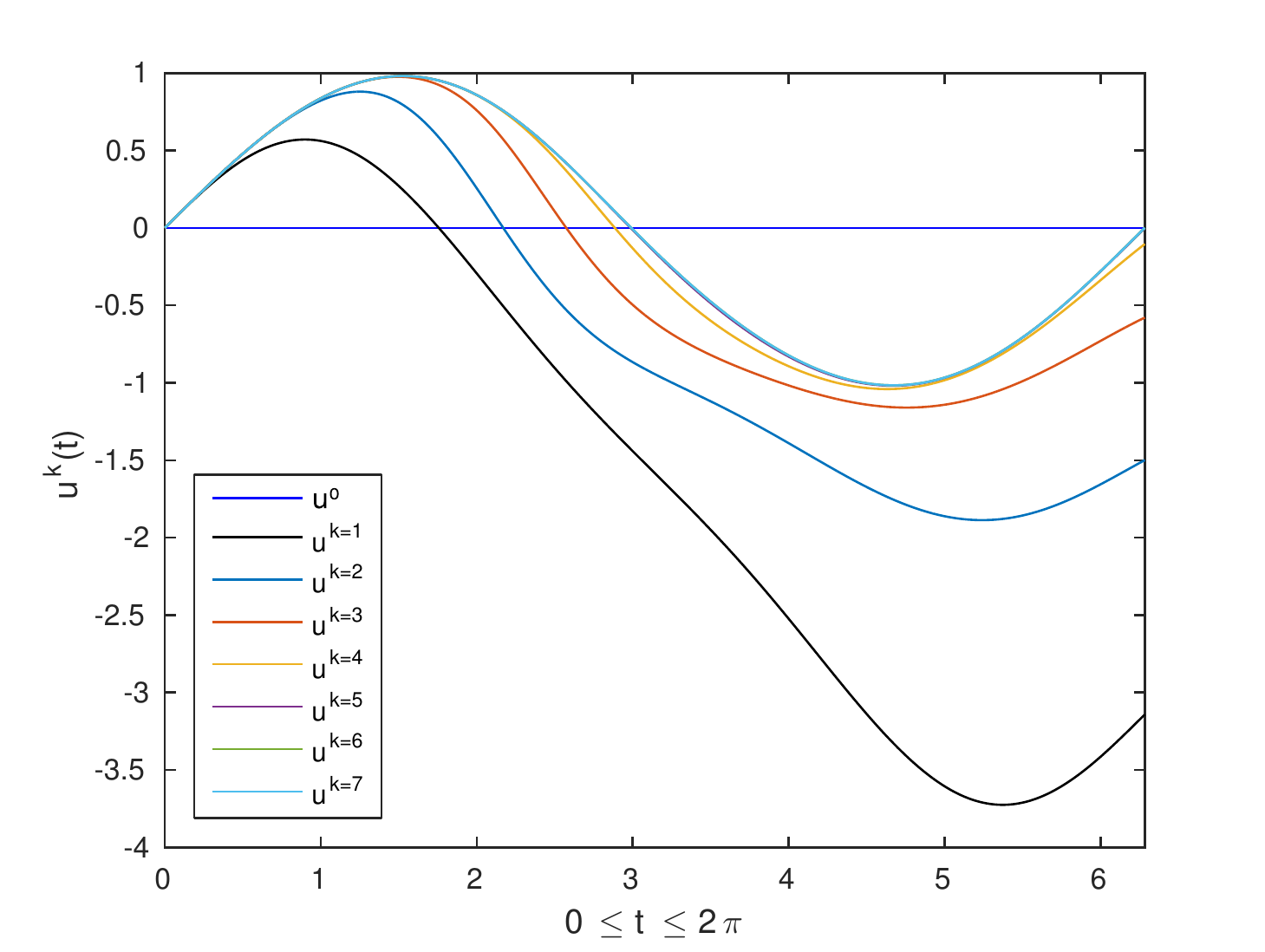} 
}%
\subfigure[First iteration]{\includegraphics[width=0.3\textwidth]{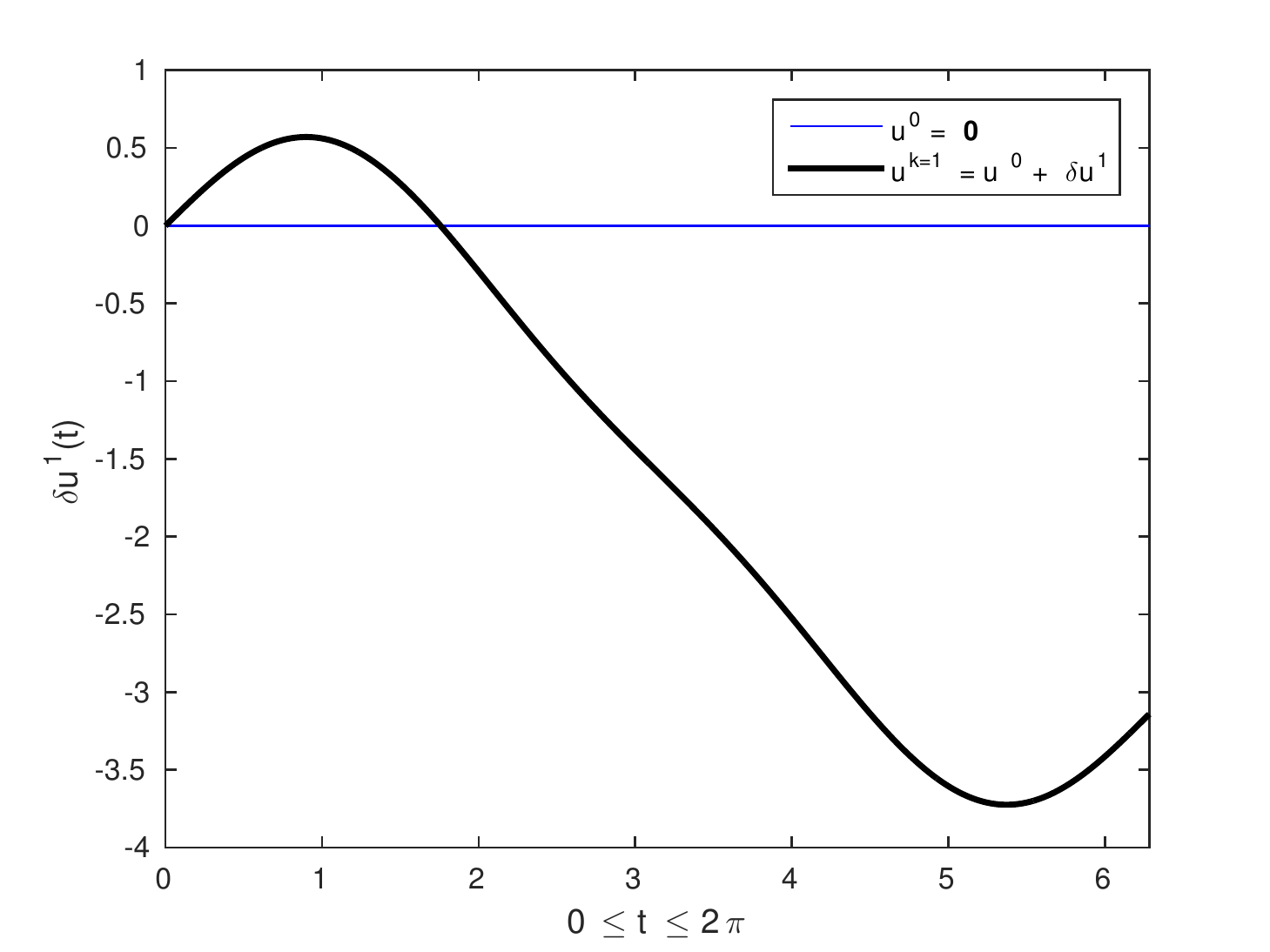}  }%
\subfigure[Additive solution]{\includegraphics[width=0.3\textwidth]{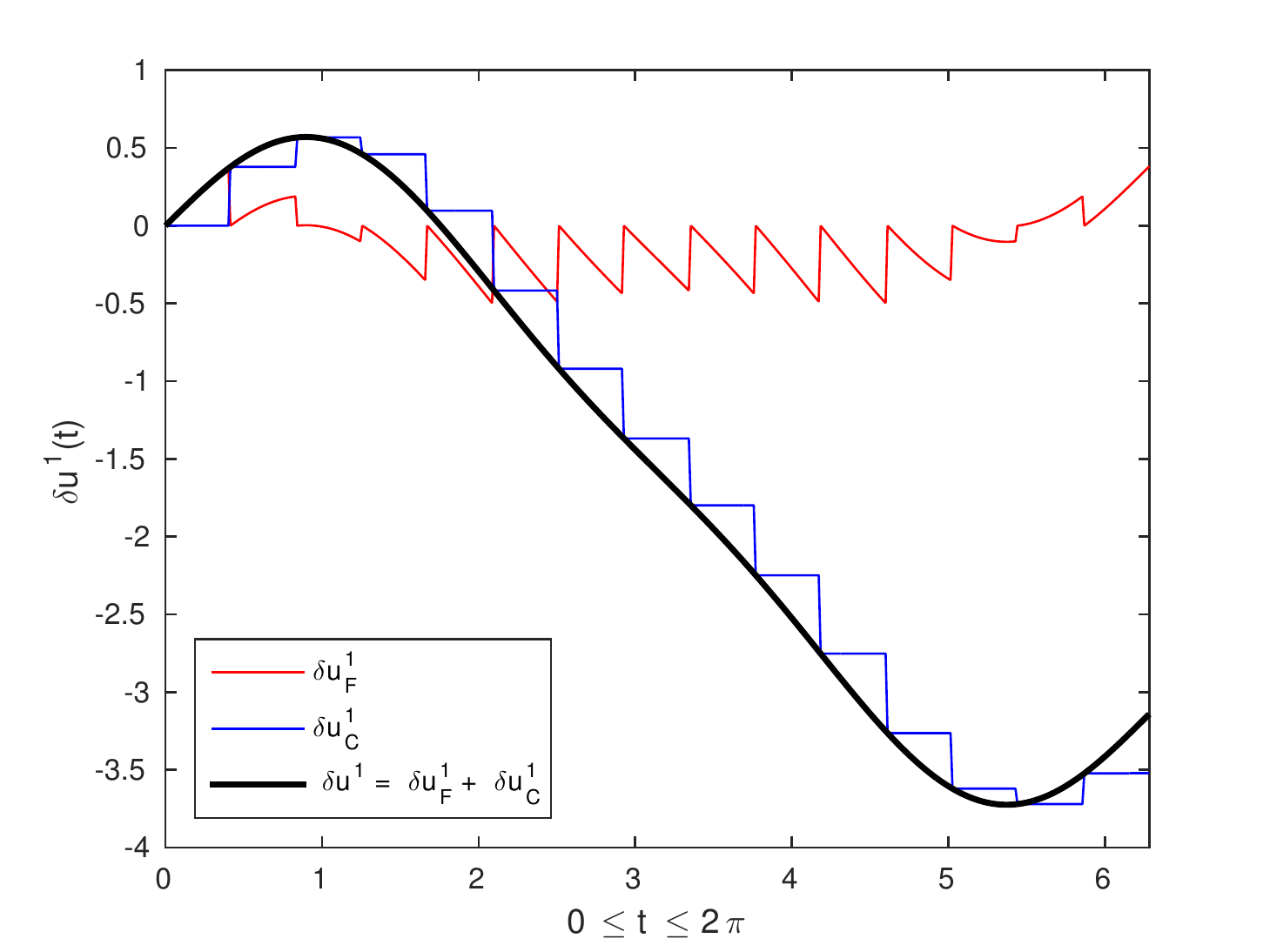}  }%
\caption{Iterations for the solution of the nonlinear equation $\partial_t u - u^2 = \rm{cos}(t) - \rm{sin}^2(t)$ on $t=[0,2\pi]$ using Newton's method and the parallel ODE direct solver. Fine and coarse solutions for the first nonlinear solution update $\delta u^1$. The time interval is discretized with 500 time steps divided into 15 subdomains. DG(0) (equivalent to Backward-Euler) is used.}

\label{fig-nonlinear_ode}

\end{figure}

{
\setalcaphskip{1ex}
\setlength{\algomargin}{1em}
\restylealgo{ruled} 
\linesnumbered
\dontprintsemicolon
\Setnlsty{textrm}{ }{:}
\SetFuncSty{\texbf}
\SetKwComment{comment}{}{}{}
\SetCommentSty{textrm}

\begin{algorithm}
  \small{
    \KwData{$\u^0$}
    \KwResult{$\u_0$ such that $\Ac_0(\u_0) = 0$}
    $\zv \leftarrow \uv^0$ \% Initial guess   \comment*[r]{$k$} 
    \While{not convergence}
    {
    Solve problem $\Jc_{0}(\zv) \yv = -\Ac_{0}(\zv)$ using the multilevel Schur complement Alg. \ref{alg:ml_schur}   with $(\Jc_{0}(\zv),-\Ac_{0}(\zv))$ \comment*[r]{$0, \ldots, \ell$}
    Assign $\zv \leftarrow \zv  + \yv$ \comment*[r]{$k$}
    }
    Return $\zv$ \comment*[r]{$k, \ldots, \ell$}  
  }
      \caption{Newton-Schur complement solver  \label{alg:newton_schur}}
  \end{algorithm}
}

\subsection{Nonlinear Schur complement-Newton methods}

Following the ideas in nonlinear domain decomposition (see, e.g., \cite{cai_nonlinearly_2002,klawonn_nonlinear_2014}), one can state the problem as a nonlinear Schur complement at a given level and next its linearization using, e.g., Newton's method.  In order to present the problem, we consider the two-level case. Later, the algorithm will be extended to multiple levels. Let us define the level-1 nonlinear Schur complement problem
$$
\Ac_1 ( \uv_1 ) = \Ac_0^\Gamma ( \Ec_0 (\uv_1) ), \quad \hbox{where } \quad  \Ec_0(\uv_1) \doteq [ \Ec_0^I(\uv_1), \uv_1 ]^T,
$$
is the nonlinear harmonic extension operator, solution of 
\begin{align}\label{nonlinear_extension}
\Ac_0^I ( \Ec_0 (\uv_1) ) \doteq \Ac_0^I ( [ \Ec_0^I(\uv_1), \uv_1 ]^T  ) = \mathbf{0}.
\end{align}
We note that the computation of $\Ec_0(\uv_1)$ involves $n_1$ independent local nonlinear ODE solvers, using the same rationale as for the linear case. We denote by $[i]$ the time steps at level-0 in the time interval $(t_1^{i},t_1^{i+1})$, i.e., $t_0^{m(i-1)+1},\ldots,t_0^{m(i)-1}$. Thus, \Eq{nonlinear_extension} can be computed as:
\begin{align}\label{seg_nl}
\Ac_0^{{[i]}} ( \Ec_0^{[i]} (\uv_1^i) ) 
 = \mathbf{0}, \quad \hbox{for} \ i =0,\ldots,n_1-1.
\end{align}
Next, we apply linearization over the level-1 nonlinear Schur complement problem, e.g., Newton's method. In this case, using the implicit function theorem (see \cite{klawonn_nonlinear_2014}), we know that
\begin{align}
  \Jc_1 ( \bar \uv_1 ) & \doteq \frac{\partial \Ac_1  }{\partial \uv_1}(\bar\uv_1) \label{jac_imp} \\
  & = \Jc_0^{\Gamma \Gamma}(\Ec_0(\bar \uv_1)) - \Jc_0^{I \Gamma}(\Ec_0(\bar \uv_1)) \Jc_0^{II}(\Ec_0(\bar \uv_1))^{-1} \Jc_0^{I \Gamma}(\Ec_0(\bar \uv_1)), \nonumber
\end{align}
i.e., the level-1 Schur complement matrix of $\Jc_0(\Ec_0(\bar \uv_1))$. The computation of the Schur complement operator can be done in parallel, as described in Sect. \ref{sec:ODEs}. As a result, the only difference between the Newton-Schur complement Alg. { \ref{alg:newton_schur}} and the nonlinear Schur complement-Newton Alg. { \ref{alg:ml_nonlinear_schur}} (for $k=1$) is the nonlinear interior correction being computed in the second case (stated in Alg. \ref{alg:ml_nonlinear_extension}). Using recursion, we can extend the nonlinear Schur complement-Newton algorithm to arbitrary levels. In Alg. \ref{alg:ml_nonlinear_schur} we state the algorithm to solve the nonlinear Schur complement at level-$k$ using Newton's method for both the nonlinear extensions in Alg. \ref{alg:ml_nonlinear_extension} and the problem itself. Thus, the resulting method involves nested Newton iterations. Let us describe these algorithms in detail.

We first consider the nonlinear harmonic extension at level-$k$. Since the values at the previous level are known, i.e., $\u_{k}^\Gamma$ is fixed, the computation of $\Ac_k(\u_k) = \mathbf{0}$ is to be computed solving local nonlinear ODE problems \Eq{seg_nl} (at level-$k$). Alg. \ref{alg:ml_nonlinear_extension} states the computation of such nonlinear ODE problem for a given level $k$ and level-$k$ time element $i$, i.e., $(t_k^i, t_k^{i+1})$. We can solve the local nonlinear ODE  using Newton's method. In any case, unless $k = 0$, we do not have an explicit expression of the nonlinear ODE. Thus, in order to compute the Jacobian (Eq. \Eq{jac_imp} for level $k$), we require to compute the $k-1$ level extension in $i$. It leads to the deployment of ${\rm card}([i])$ nonlinear ODEs at the level $k-1$ in elements $j \in [i]$. Again, if the next level is not the level-0, we have to proceed recursively to solve these local problems. All these tasks are described in Alg. \ref{alg:ml_nonlinear_extension}.
{
\setalcaphskip{1ex}
\setlength{\algomargin}{1em}
\restylealgo{ruled} 
\linesnumbered
\dontprintsemicolon
\Setnlsty{textrm}{ }{:}
\SetFuncSty{\texbf}
\SetKwComment{comment}{}{}{}
\SetCommentSty{textrm}

\begin{algorithm}
  \small{
    \KwData{$k$, $i$, $\v$, $\vv_0$}
    \KwResult{$\Ec^{[i]}_{k}(\v)$, $\Jc_k^{[i]}(\Ec_k^{[i]}(\v))$, $-\Ac_k^{[i]}(\Ec_k^{[i]}(\v))$}
    \uIf{ $k = 0$ }
    {
      Solve $\Ac_0^{[i]}(\Ec_0^{[i]}(\v)) = \mathbf{0}$ using Newton's method (local nonlinear ODE), in order to get $\Ec_0^{[i]}(\v)$ \comment*[r]{0}
      Compute $\Jc_0^{[i]}(\Ec_0^{[i]}(\v))$, $-\Ac_0^{[i]}(\Ec_0^{[i]}(\v))$ \comment*[r]{0}
    }
    \Else
    {
      \% Solve $\Ac_{k}^{[i]}(\Ec_{k}^{[i]}(\v))) = \mathbf{0}$ using Newton's method as follows \;
      $\zv^{[i]} = \vv_0^{[i]}$ \% Initial guess   \comment*[r]{$k$} 
      \While{not convergence}
      {
        \For{$j \in [i]$}
        {   
          Solve
          $\Ac_{k-1}^{[j]}( \Ec_{k-1}^{[j]}(\zv^{j})
          ) = \mathbf{0}$
          using Alg. \ref{alg:ml_nonlinear_extension} with
          $(k-1,j,\zv^{j},\vv_0)$ \comment*[r]{$k-1$} 
          Compute $\Jc_{k-1}^{[j]}( \Ec_{k-1}^{[j]}(\zv^j))$ and $- \Ac_{k-1}^{[j]}(\Ec^{[j]}_{k-1}( \zv^j ))$ \nllabel{lin-local-value} \comment*[r]{$k-1$} 
        }
        Assemble  $\Jc_k^{[i]}(\zv )$ and $- \Ac_k^{[i]}( \zv )$ using the local contributions in l.\ref{lin-local-value} and formula \Eq{jac_imp} \comment*[r]{$k-1\to k$} 
        Solve $\Jc_k^{[i]}( \zv ) 
         \yv = - \Ac_k^{[i]}( \zv )
        $ and assign $\zv \leftarrow \zv  + \yv$ \comment*[r]{$k$}  
      }
      Return $( \zv, \ \Jc_k^{[i]}(\zv ), - \Ac_k^{[i]}( \zv ))$ \comment*[r]{$k$} 
    }
      \caption{Nonlinear harmonic extension \label{alg:ml_nonlinear_extension}}}
  \end{algorithm}
}

Once we have defined the nonlinear extension operator, we can state the level-$k$ nonlinear Schur complement solved with Newton's method. Using the expression of the Jacobian in \Eq{jac_imp} (for level $k$), we require to compute the nonlinear extension at the next level, using Alg. \ref{alg:ml_nonlinear_extension} at all level-$k$ elements. It leads to a level-$k$ linear ODE to be solved, for which we can use the multilevel Schur complement approach in Alg. \ref{alg:ml_schur}, exploiting levels $k, \ldots, \ell$ in its solution. As commented above, the level-1 nonlinear Schur complement consists in Newton iterations like in Alg. \Eq{alg:newton_schur}, but with the main difference that one performs a nonlinear interior correction of the values at level-0, solving the nonlinear ODE locally. The level-$k$ nonlinear Schur complement requires nested nonlinear iterations to perform the nonlinear harmonic extension at that level, but the Jacobian problem to be solved only involves levels $k,\ldots,\ell$. As commented above, other linearization techniques can readily be used.
{
\setalcaphskip{1ex}
\setlength{\algomargin}{1em}
\restylealgo{ruled} 
\linesnumbered
\dontprintsemicolon
\Setnlsty{textrm}{ }{:}
\SetFuncSty{\texbf}
\SetKwComment{comment}{}{}{}
\SetCommentSty{textrm}

\begin{algorithm}
  \small{
    \KwData{$k$, $\u^0$, $\uv^0_0$}
    \KwResult{$\uv_0$ such that $\Ac_0(\uv_0) = 0$}
    $\zv_0 = \uv_0^0$ \% Initial guess   \comment*[r]{$k$} 
    \While{not convergence}
    {
      \For{$i = 0, \ldots, n_{k}$}
      {
        Compute nonlinear harmonic extension $\Ec_{k-1}^{[i]}(\zv^{i}_k)$ and the local contributions
        $\Jc_{k-1}^{[i]}(\Ec_{k-1}^{[i]}(\zv^i_k))$, $-\Ac_{k-1}^{[i]}(\Ec_{k-1}^{[i]}(\zv^i_k))$ using Alg.~\ref{alg:ml_nonlinear_extension} with $(k-1,i,\zv^i_k,\zv)$ \nllabel{lin-local-value2}  \comment*[r]{$k-1, \ldots, 0$} 
      }
    Assemble the level-$k$ Schur complement around $\Ec_{k-1}(\zv_k)$ using the local contributions in l.\ref{lin-local-value2} and formula \Eq{jac_imp}  \comment*[r]{$k-1\to k$} 
    Solve problem $\Jc_{k}(\zv_k) \yv = -\Ac_{k}(\zv_k)$ using the multilevel Schur complement Alg. \ref{alg:ml_schur}   with $(\Jc_{k}(\zv_k),-\Ac_{k}(\zv_k))$ \comment*[r]{$k, \ldots, \ell$}
    Assign $\zv_k \leftarrow \zv_k  + \yv$ \comment*[r]{$k$}
    }
    Return $\zv_0 $ \comment*[r]{$k, \ldots, \ell$}  
  }
      \caption{level-$k$ nonlinear Schur complement-Newton solver  \label{alg:ml_nonlinear_schur}}
  \end{algorithm}
}

\section{Numerical experiments}\label{sec:numerical_experiments}

\subsection{Experimental set-up}
In this section we evaluate the weak scalability of the proposed
methods. We consider the time interval $(0,T]$, which is divided into
$n_0$ time elements. The parallel solver relies on a level-1 coarser
time partition into $n_1$ time elements (time subdomains); thus, every
level-1 time element is defined by aggregation of $\frac{n_0}{n_1}$
level-0 time elements. $n_1$ processors are used for the simulations
in all cases. One higher level of coarsening gives a partition into
$n_2$ elements; analogously, every level-2 time element is defined by
aggregation of $\frac{n_1}{n_2}$ level-1 time elements. See
Fig. \ref{fig-part} for a graphical illustration. Only a subset of
$n_2$ processors have duties at level-2. Two-level methods involve
level-0 and level-1 duties, whereas the three-level methods also
involve level-2 duties. Higher levels have not been needed at the
scales considered below but could be needed at largest scales to keep
good weak scalability.

Three different approaches are considered for solving the
{ Lotka}-Volterra system of nonlinear ODEs: 1) the pure sequential
approach (no parallelism is exploited, using linearization at every
time step), 2) the Newton-Schur complement approach (that involves a
global linearization) in Alg. \ref{alg:newton_schur}, and 3) the
1-level nonlinear Schur complement-Newton approach in
Alg. \ref{alg:ml_nonlinear_schur}. Time integration is performed with
the Backward Euler scheme with a constant time step. The linearization
of the nonlinear problems is performed with a hybrid Picard-Newton
technique, where Newton's method is activated when the discrete
$L_2$-norm of the residual is below $10^2$. (Even though all the
algorithms have been stated using the full Newton's method for
simplicity, its extension to Picard linearization (and hybrid
approaches) is straightforward.)  Sequentiality is exploited at the
local level in each processor for all the different approaches
(level-0 and higher level systems are lower block-triangular; see
Eqs. \ref{global_ODE} and \ref{interface_correction},
respectively). The stopping criteria for the sequential and
Newton-Schur complement solvers is the reduction of the discrete
$L_2$-norm of the nonlinear local/global residual below an absolute
value of $10^{-8}$. For the third approach, the global residual
tolerance is fixed to $10^{-8}$, while local nonlinear problems and
Schur complement solutions converge below $10^{-10}$.

The Lotka-Volterra equations are a system of nonlinear ODEs frequently
used to describe the dynamics of biological systems in which two
species interact, one as a predator and the other one as a prey. Our
unknown functions, namely $(u,v)$, represent the evolution in time of
the number of units of each one of the species considered. The system
of nonlinear, first order, differential equations reads
\begin{align}\label{eq-lotka_volterra}
  \displaystyle\frac{{\rm d}u(t)}{{\rm d}t} = \alpha u - \beta u v  \\
  \displaystyle\frac{{\rm d}v(t)}{{\rm d}t} = \delta uv - \gamma v 
\end{align}
where $\alpha,\beta,\gamma,\delta$ are positive parameters that model the interaction of the two species. The system is solved for $t\in(0,T]$. Appropriate initial conditions (initial number of preys and predators) must be provided, i.e., $u(0), v(0)$. Out of these equations, and graphically illustrated with an example in Fig. \ref{fig-lotka_volt_example}, if no interaction between the species is modelled (uncoupled linear equations), the preys grow in number while the predators decrease \ref{fig-no_interaction}. In the second case, the species interaction leads to the frequency plot \ref{fig-interaction}. 
 
 \begin{figure}[h!!]
\centering
\subfigure[No interaction]{\includegraphics[width=0.4\textwidth]{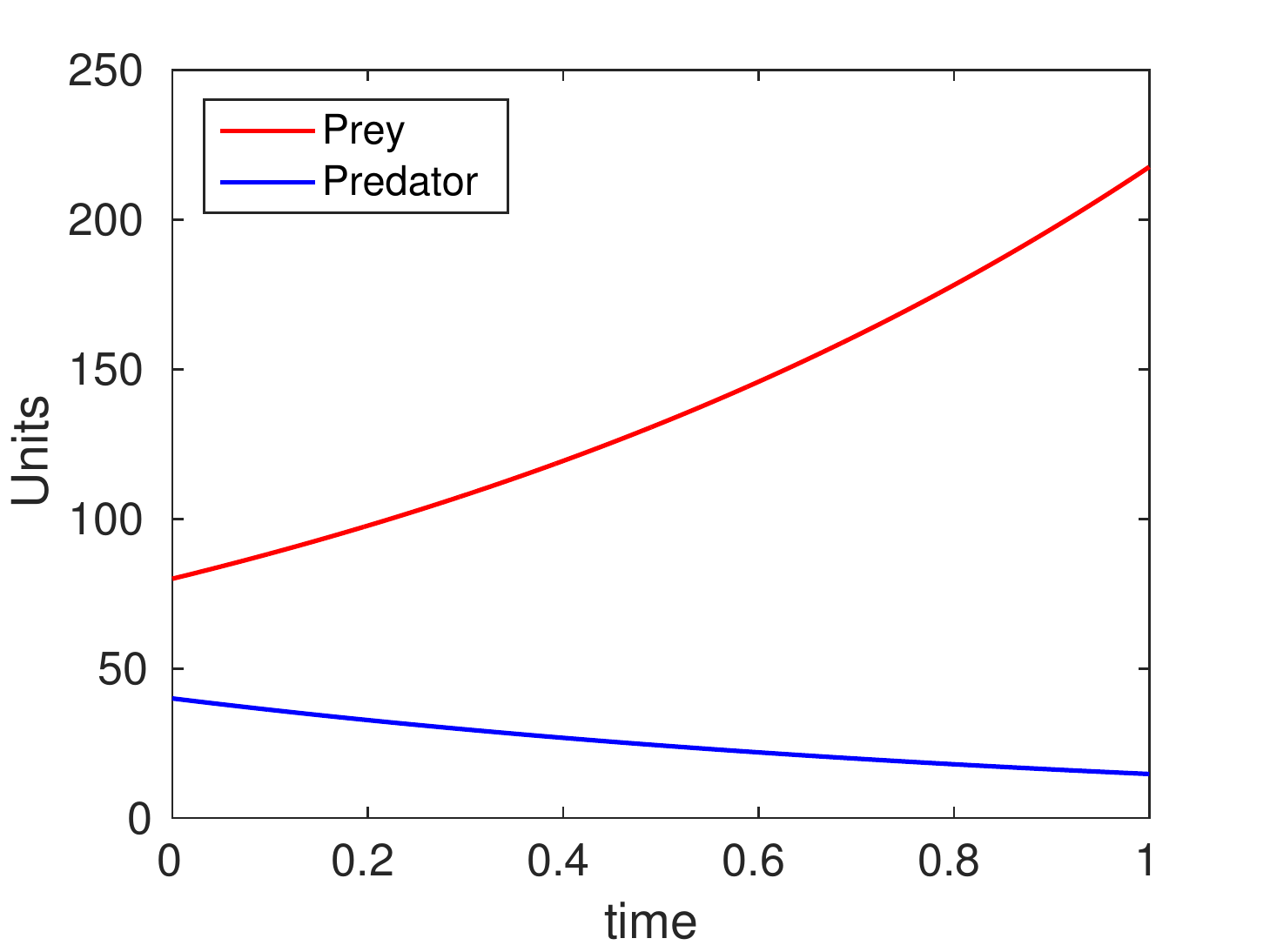} 
\label{fig-no_interaction}}
\subfigure[$\alpha=3$, $\beta=0.2$, $\delta=0.1$, $\gamma=2$]{\includegraphics[width=0.4\textwidth]{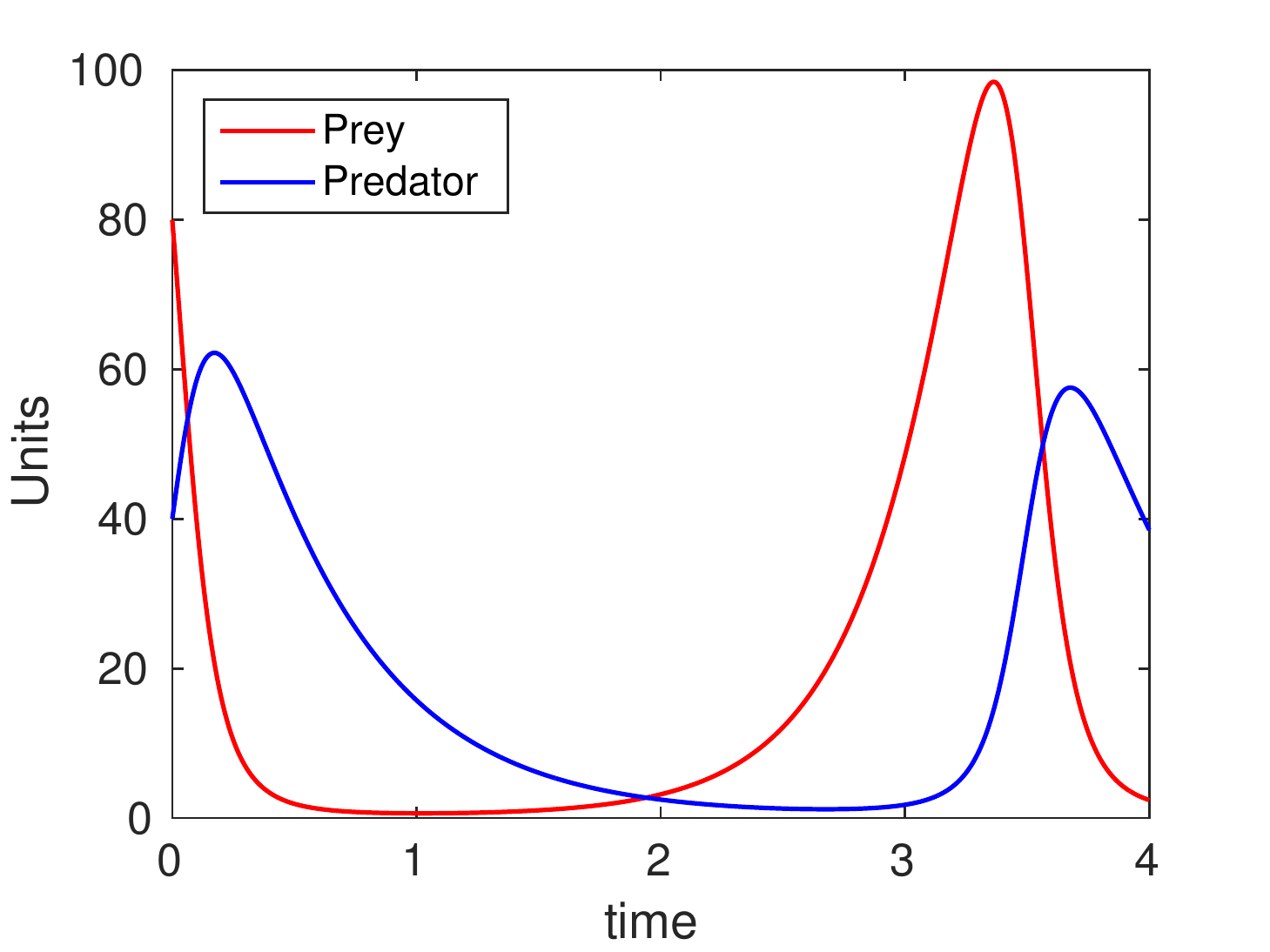}\label{fig-interaction}  }%
\caption{Solution to the Lotka-Volterra equations when different scenarios are considered. } \label{fig-lotka_volt_example}
\end{figure} 

\subsection{Two-level solvers}
In this subsection, weak scalability results are presented for the
solution of the Lotka-Volterra equations with a two-level method. The
local problem size at level-0 is fixed to $\theta =
\frac{n_0}{n_1}$. We consider a weak scalability analysis in which we
keep fixed the local problem size $\frac{n_0}{n_1}$ and increase
$n_1$, i.e., the number of subdomains (level-1 time elements). Thus,
we are increasing the global time steps being used to solve the
ODE. The Schur complement at level-1 has size $n_1$. We stop the
analysis when $n_1 \approx \frac{n_0}{n_1} $, since we would require a
three-level algorithm (see Sect. \Eq{sec:par_ef}). 
 
Figs. \ref{fig-Ex50} and \ref{fig-Ex100} show a comparison between the weak scalability of the sequential solver (1 processor performs sequentially the full computation) and the parallel solvers for three different local problem sizes $\frac{n_0}{n_1}$. In these plots, the parallel solvers computing times are an aggregation of the local solvers (level-0) and the coarse solver (level-1) CPU times. The number of iterations for the sequential approach is an average value for all nonlinear time step problems. For parallel approaches, it shows the number of global nonlinear iterations required to meet the convergence requirements.   

 \begin{figure}[h!!]
\centering
\subfigure[Newton-Schur complement solver]{\includegraphics[width=0.32\textwidth]{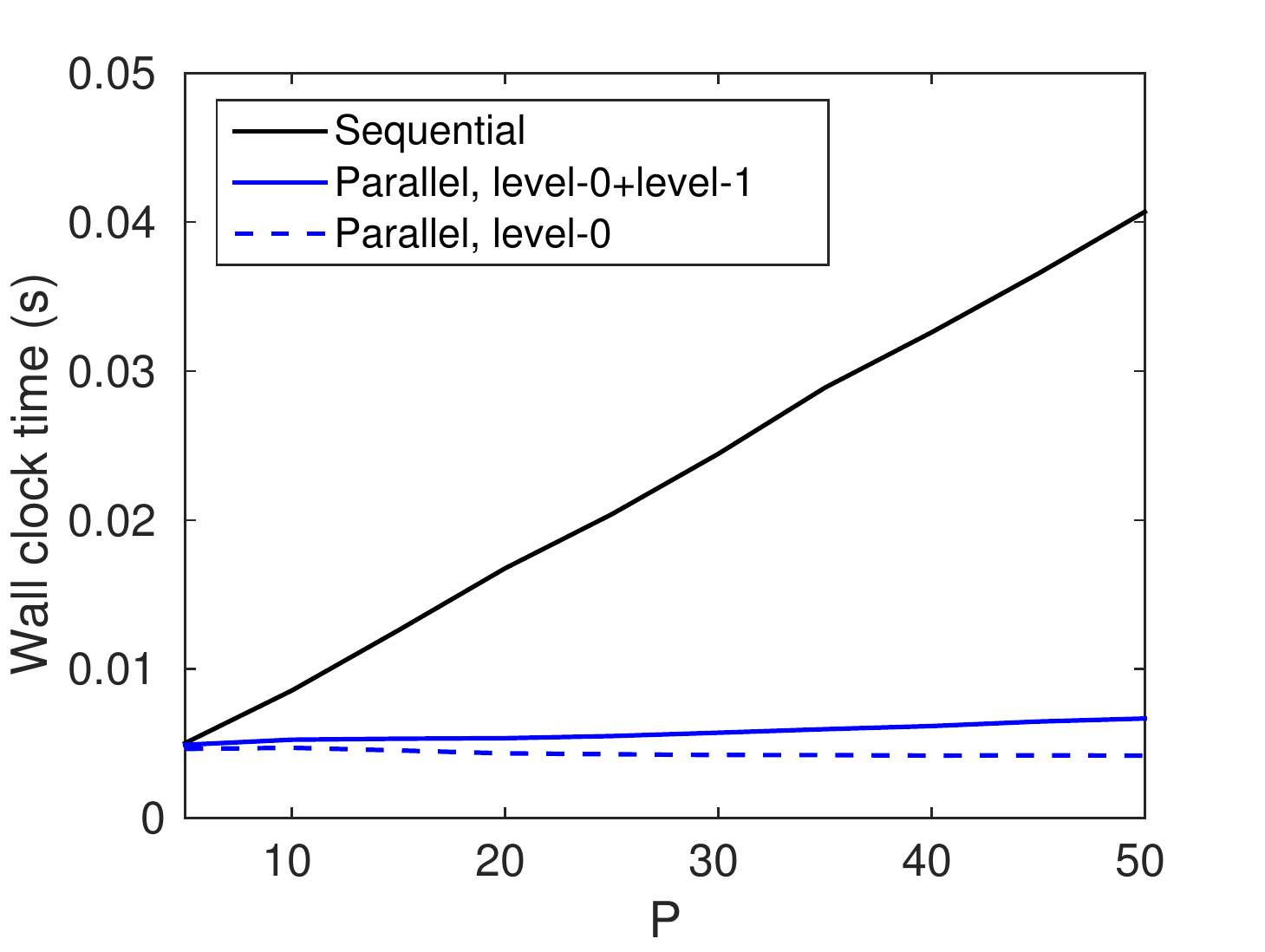} 
\label{fig-nlg_50}}
\subfigure[Nonlinear Schur complement solver]{\includegraphics[width=0.32\textwidth]{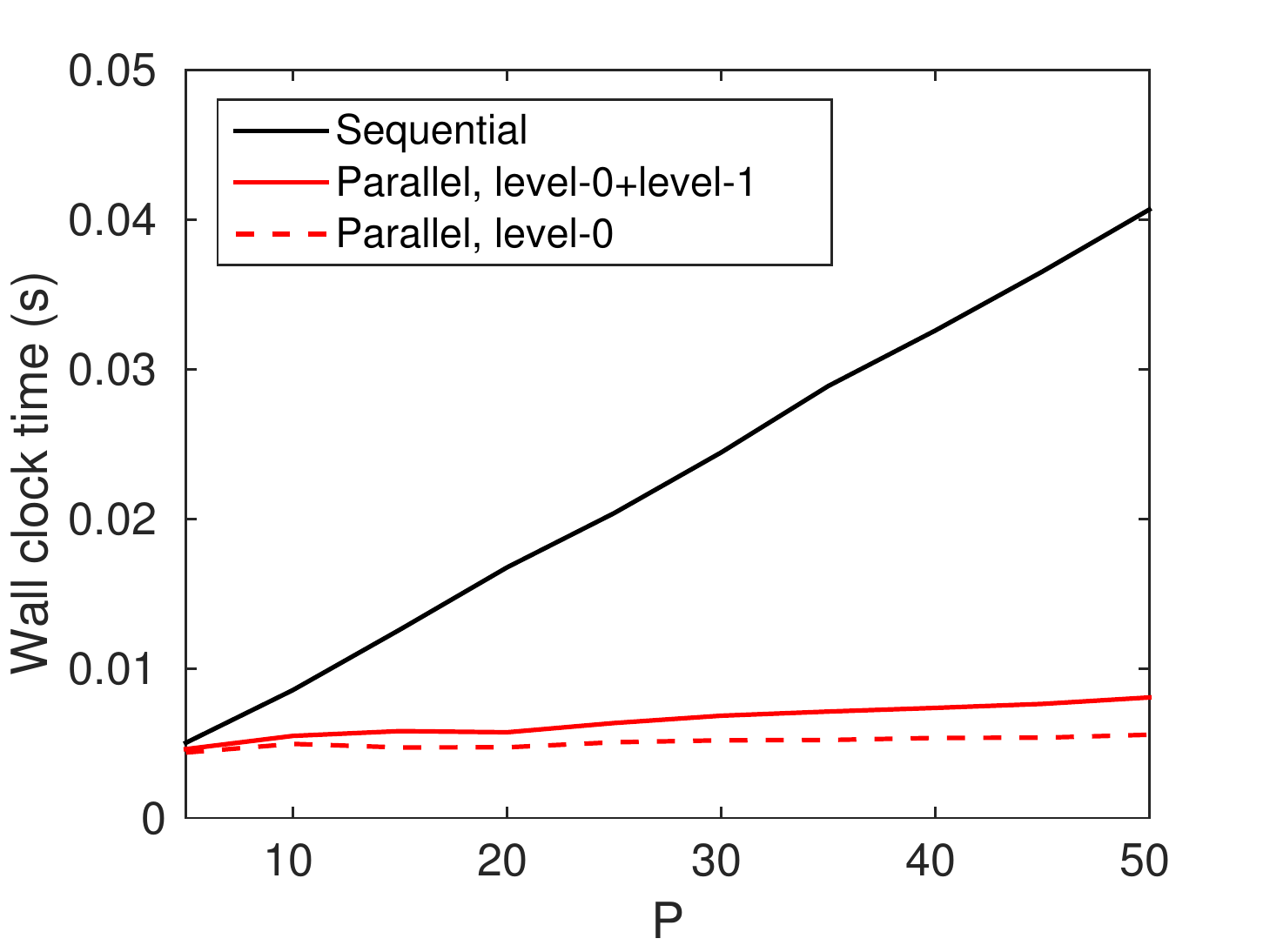}\label{fig-nls_50}   }%
\subfigure[Number of iterations]{\includegraphics[width=0.32\textwidth]{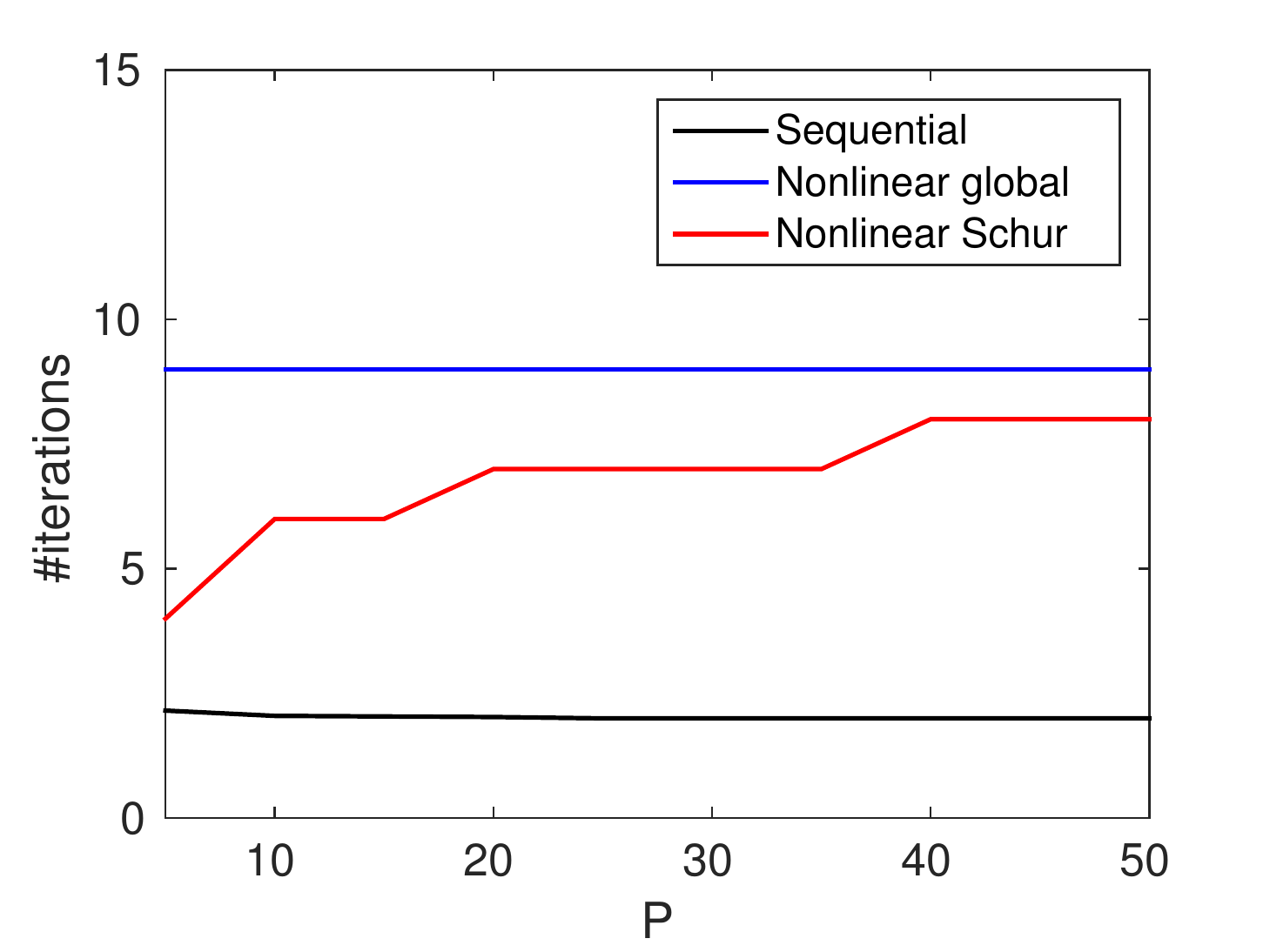}\label{fig-its_50}   }%
\caption{Weak scalability for local problem size $n_0/n_1=50$. The problem is solved with $\alpha=3$, $\beta=0.2$, $\delta=0.1$, $\gamma=2$ and $t\in(0,3]$. Initial guess $u(0)=10$, $v(0)=40$.  } 
\label{fig-Ex50}
\end{figure} 

 \begin{figure}[h!!]
\centering
\subfigure[Newton-Schur complement solver]{\includegraphics[width=0.32\textwidth]{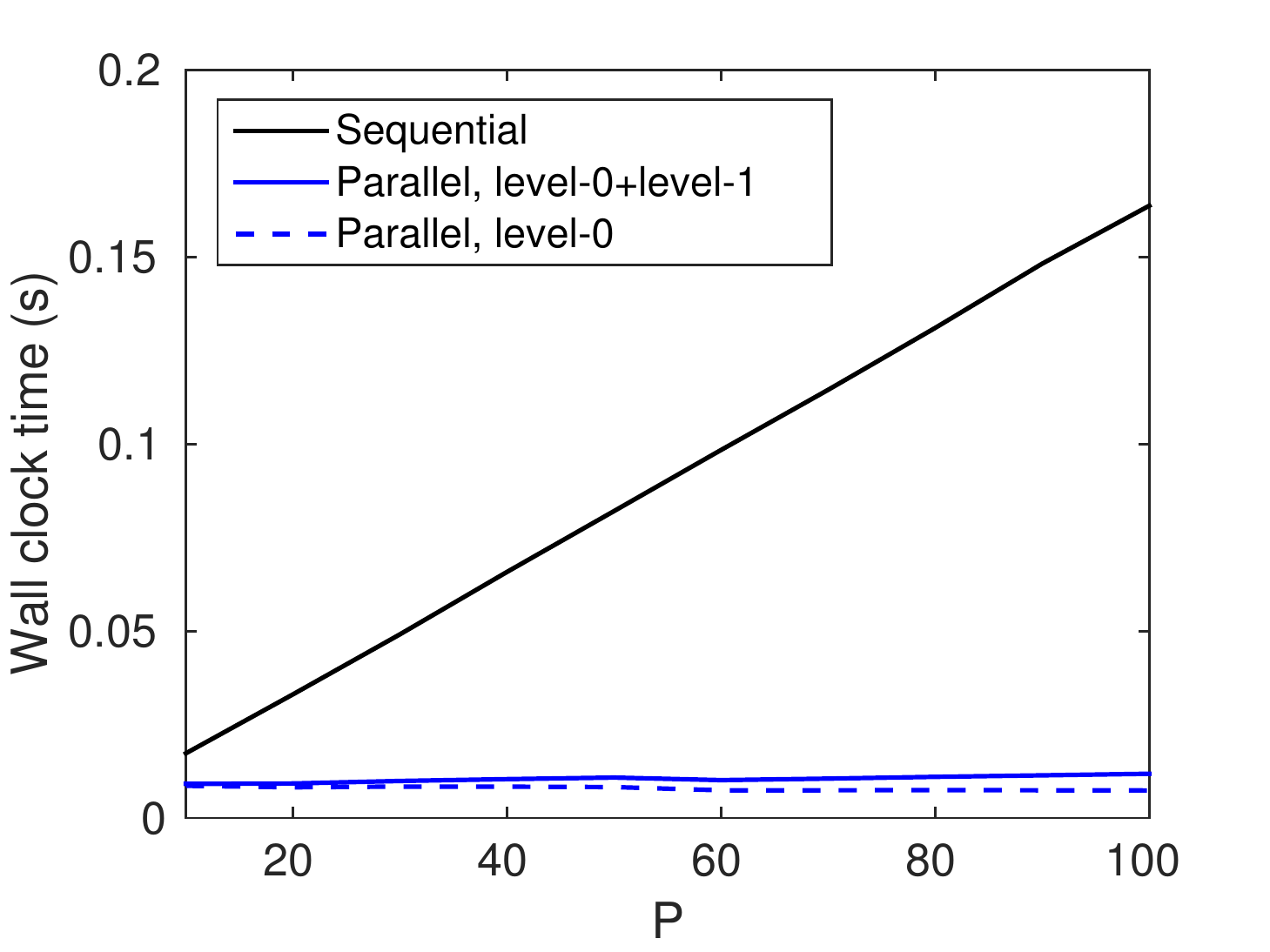} 
\label{fig-nlg_100}}
\subfigure[Nonlinear Schur complement solver]{\includegraphics[width=0.32\textwidth]{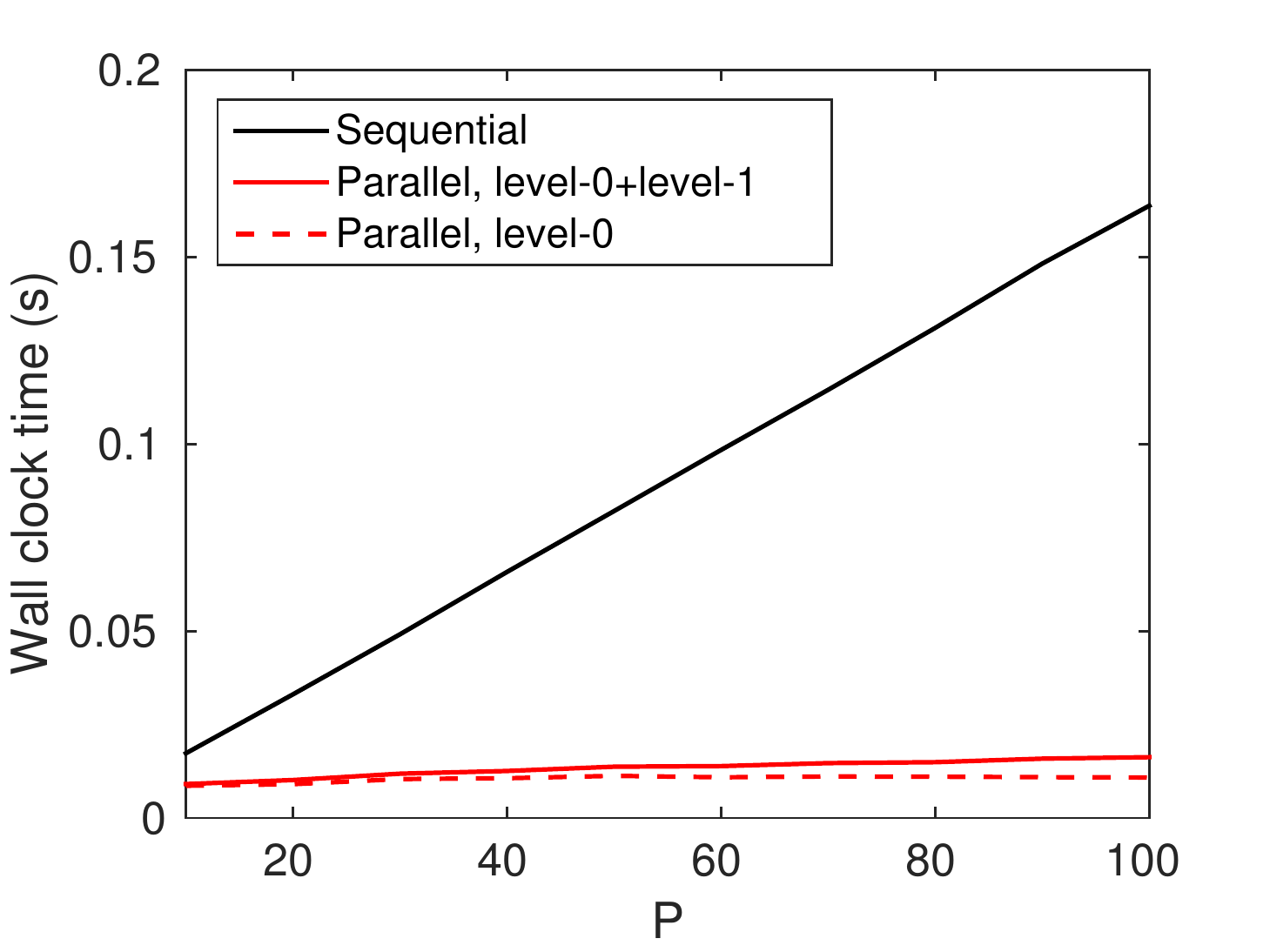}\label{fig-nls_100}   }%
\subfigure[Number of iterations]{\includegraphics[width=0.32\textwidth]{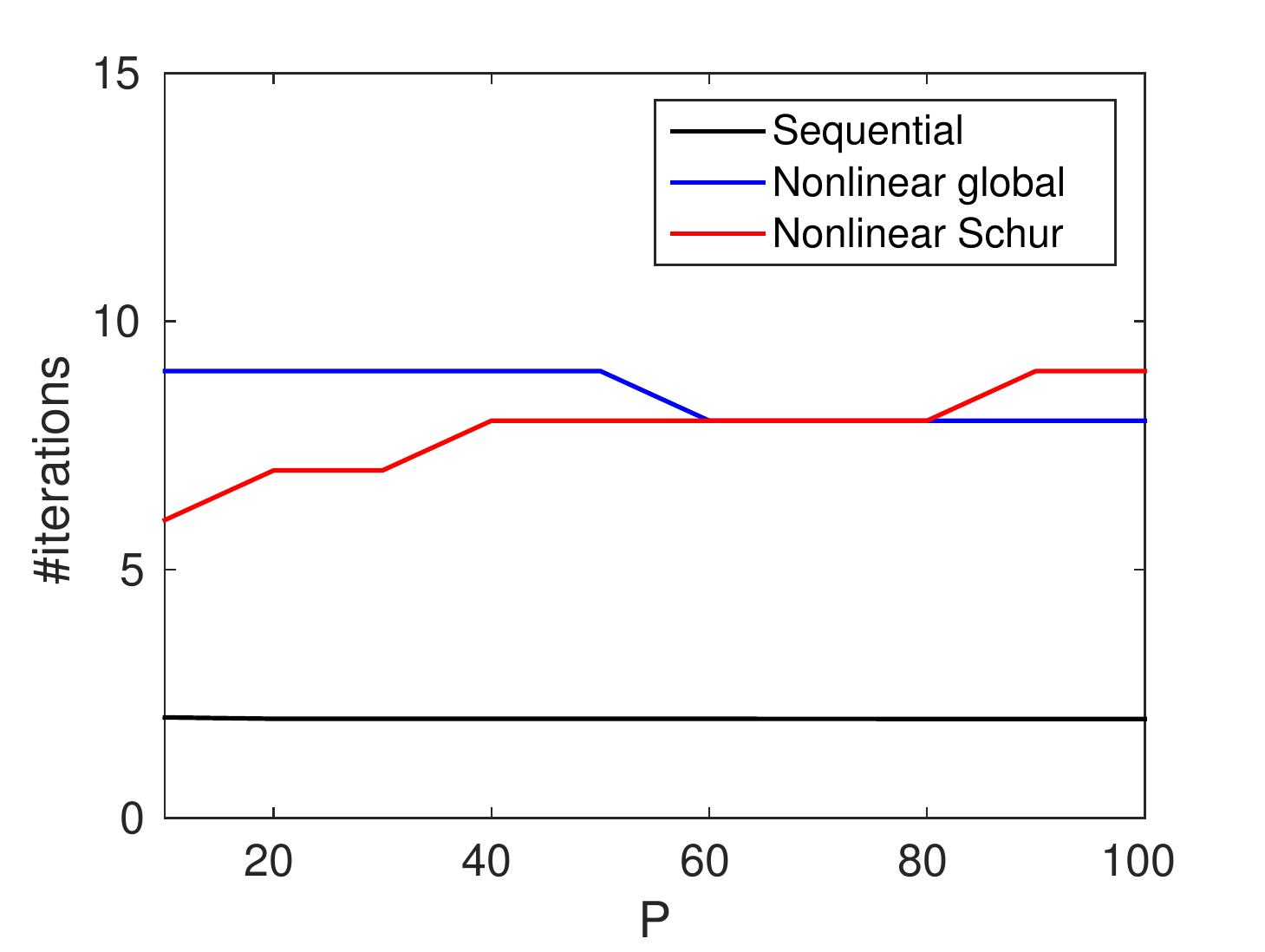}\label{fig-its_100}   }%
\caption{Weak scalability for $n_0/n_1=100$. The problem is solved with $\alpha=3$, $\beta=0.2$, $\delta=0.1$, $\gamma=2$ and $t\in(0,3]$. Initial guess $u(0)=10$, $v(0)=40$  } 
\label{fig-Ex100}
\end{figure} 


Out of these plots, we can draw some conclusions. First, parallel approaches reduce from the very beginning the time-to-solution of the simulations. Second, and most important, excellent weak scalability is observed for both parallel solvers; we can solve X times more time steps increasing X times the number of tasks, in approximately the same CPU time.

Regarding the Newton-Schur complement approach, one can observe that the CPU time at level-1 is always less than half the CPU time at level-0. It is due to the fact that we stop the scalability analysis before the coarse space is larger than the level-0 local problems and the fact that two local solvers are required at level-0 (see Alg. \ref{alg:ass_schur}) per nonlinear iteration. For the nonlinear Schur complement-Newton approach, the number of global iterations to converge is lower in most cases, but it involves nested nonlinear iterations. The Newton-Schur complement approach is slightly faster than the nonlinear Schur complement-Newton one. In any case, it can strongly be affected by the tolerances being chosen in the nonlinear Schur complement-Newton method. Clearly, increasing the load per processor (local problem sizes), the benefit of the parallel solvers compared to the sequential approach becomes more notorious.

Results are presented up to the limit of $n_1\approx \frac{n_1}{n_0}$, where the coarse problem size grows above the local problem sizes. Aiming to exploit further concurrency, the next section is devoted to show numerical results in a multilevel approach.

\subsection{Multilevel solver}

The multilevel (three-level) approach is be activated when the coarse problem size exceeds the local problem size. The first coarsening leads to a computation with $n_1$ level-0 local problems of size $\frac{n_0}{n_1}$ and a coarse problem of size $n_1$ (level-1 partition), such that $n_1 > \frac{n_0}{n_1}$. At this point, another level of coarsening is introduced to exploit further concurrency. We consider a coarsening of the level-1 time partition into $n_2$ local problems of size $\frac{n_1}{n_2}$ and a coarse problem of size $n_2$ such that $n_2 < \frac{n_1}{n_2}$ (beyond this limit, a four-level method would be required). In the following plots, the parallel solvers computing times are an aggregation of the level-0 local solver CPU time ($n_1$ parallel tasks), the level-1 local solver CPU time ($n_2$ parallel tasks) and the level-2 global solver CPU times (Schur complement sequential solve in one processor). 

In Fig. \ref{fig-l3}, CPU times for the three-level (i.e., $\ell=2$) Schur complement solve are shown. Out of the plots, it can be observed that local solves computing times for the level-0 and level-1 tasks are of the same order, since local problem sizes is kept constant for both levels. The Schur complement computation (level-2) time grows with the number of processors, as expected. Again, the level-2 the Schur complement CPU time  in Fig. \ref{fig-l3} does not exceed half the CPU time of level-0 and level-1 local solves CPU time, due to the same reasons commented above, since  $n_2 < \frac{n_1}{n_2}$.  Beyond this limit, i.e.,  for $n_2 > \frac{n_1}{n_2}$, another level would be needed.
 \begin{figure}[h!!]
\centering
\subfigure[$n_0/n_1=n_1/n_2=10$]{\includegraphics[width=0.32\textwidth]{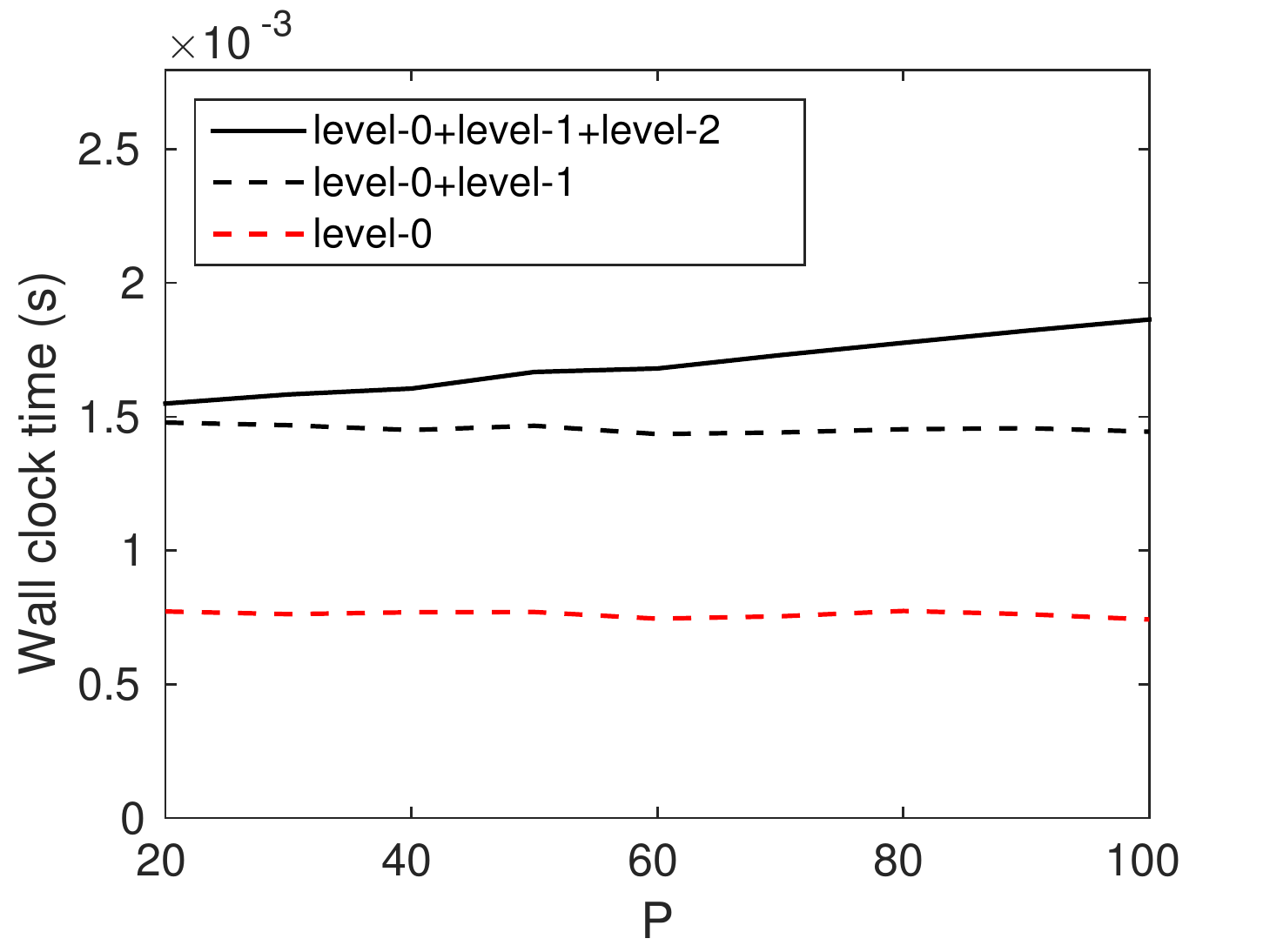} 
\label{fig-l3_10}}
\subfigure[$n_0/n_1=n_1/n_2=20$]{\includegraphics[width=0.32\textwidth]{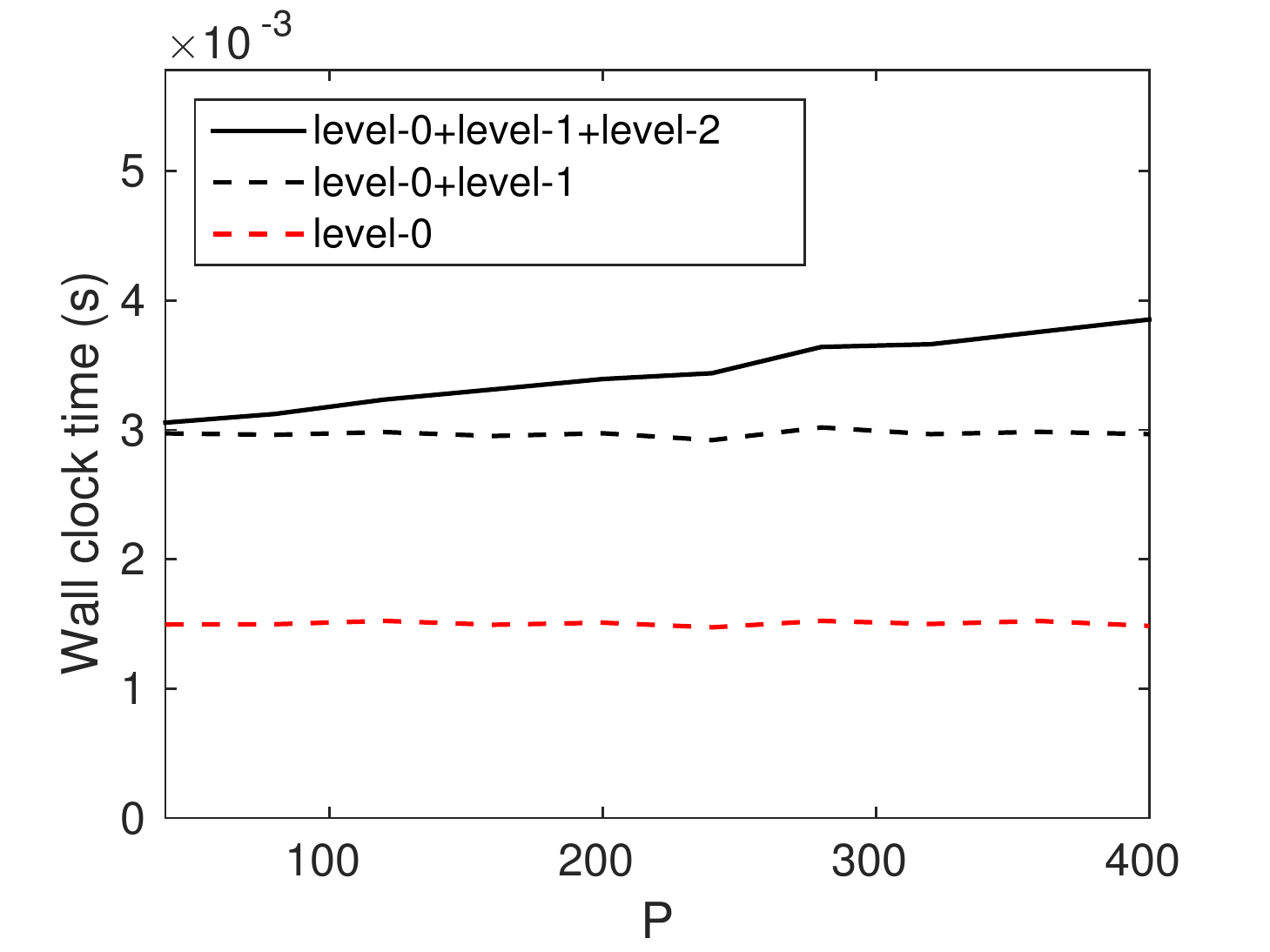}\label{fig-l3_20}   }%
\subfigure[$n_0/n_1=n_1/n_2=50$]{\includegraphics[width=0.32\textwidth]{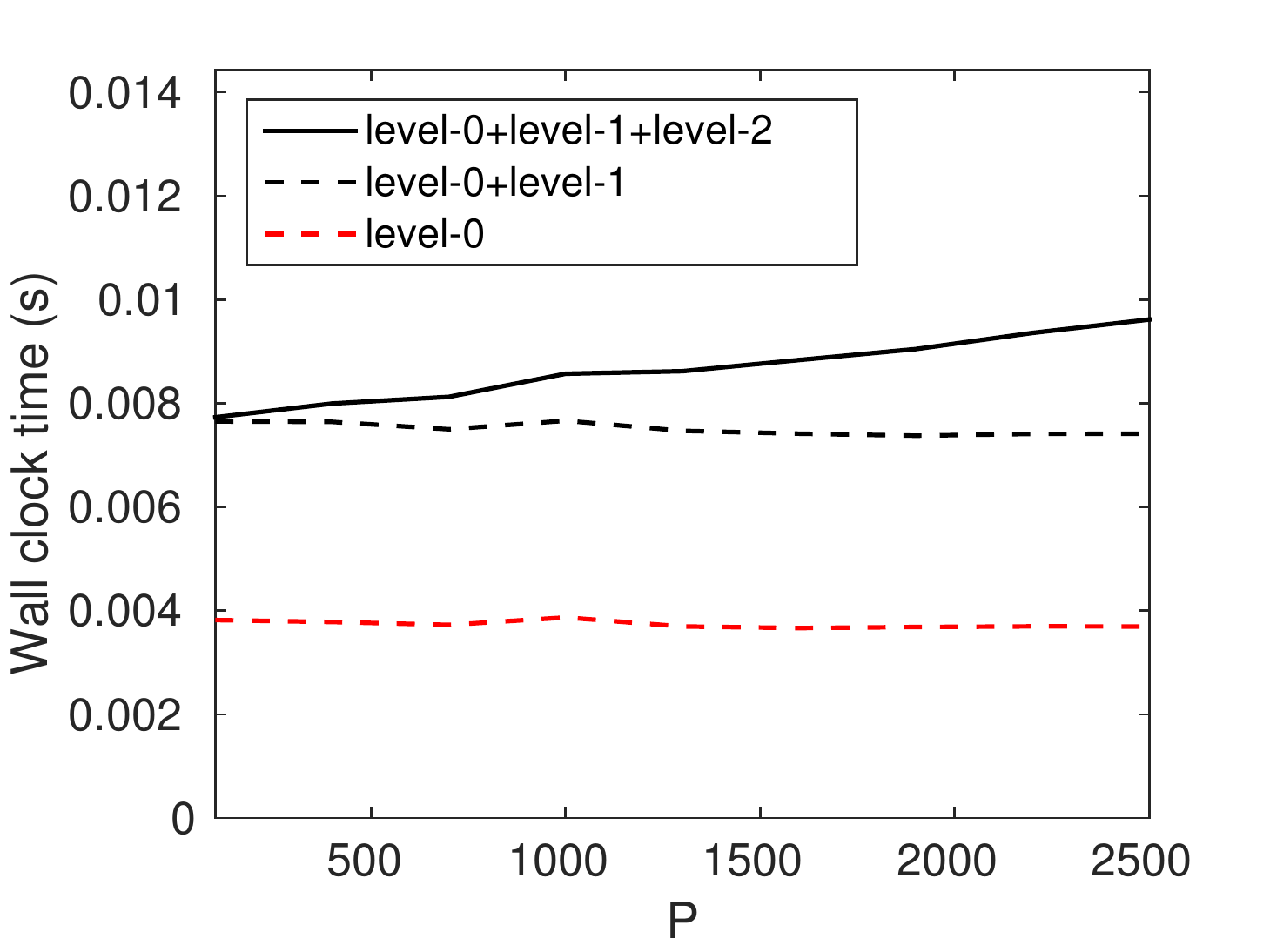}\label{fig-l3_50}   }%
\caption{Three-level ($\ell=2$) Schur complement solver with different local sizes.} 
\label{fig-l3}
\end{figure} 

In Fig. \ref{fig-l2_l3}, a comparison between the two-level ($\ell=1$) and the three-level ($\ell=2$) approaches is presented. In the first part of the plot, results are shown for the two-level approach only, since the size of the level-1 problem is still below the local sizes of the level-0 problems. The three-level technique is activated when the size of the coarse problem is two times the size of the local problems of the finest level, leading to a level-2 partition into two time elements. Out of the plots, the multilevel approach shows much better efficiency than the two-level approach. It is important to note that to include additional levels do not affect nonlinear convergence of Newton-Schur and 1-level nonlinear Schur complement-Newton methods, but it has benefits in the computation of the linear ODE systems. As a result, more levels show better computing times in these approaches.

 \begin{figure}[h!!]
\centering
\subfigure[$n_0/n_1=50$, $n_1/n_2=50$]{\includegraphics[width=0.45\textwidth]{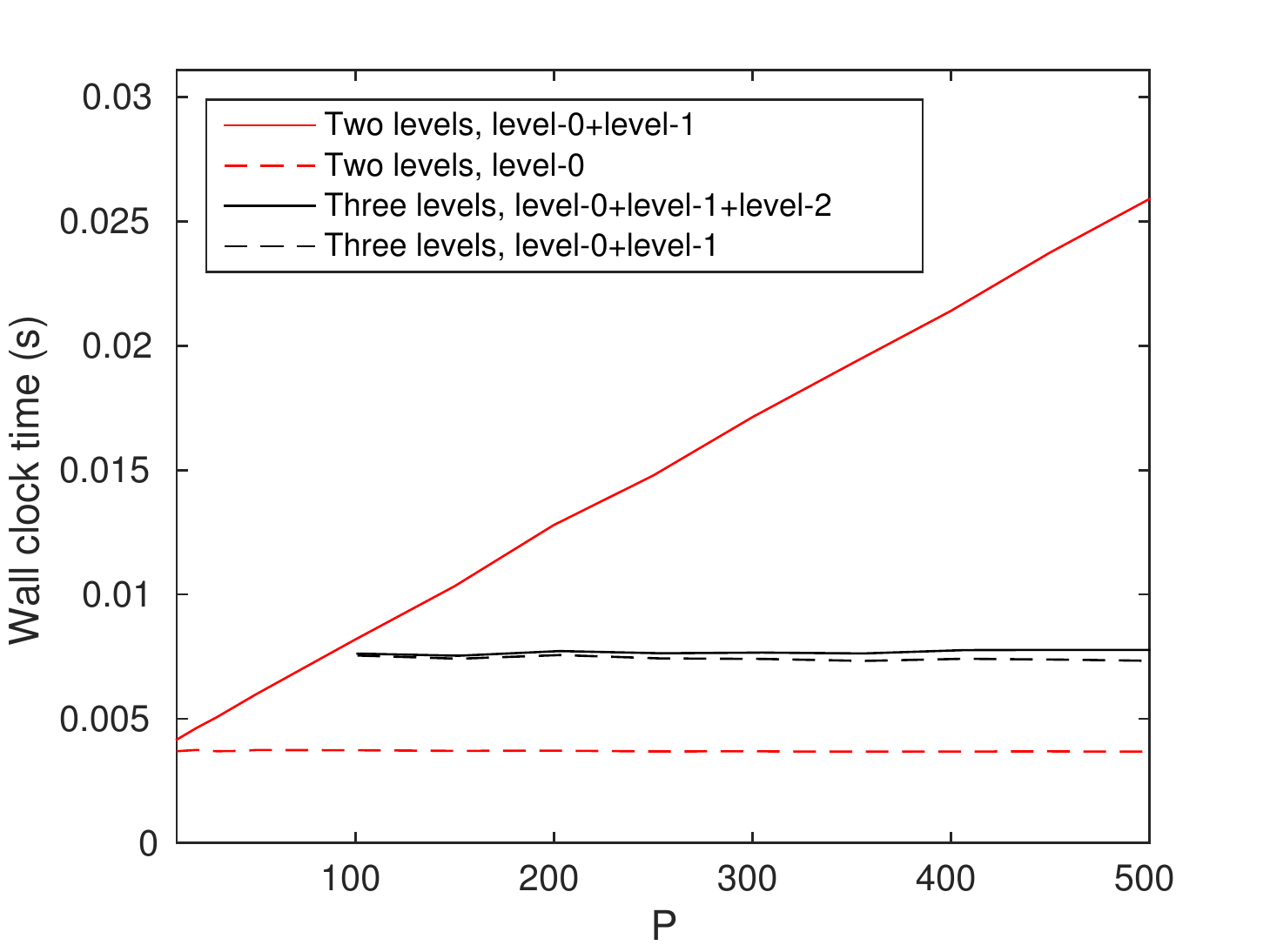} 
\label{fig-l2_l2_50}}
\subfigure[$n_0/n_1=100$, $n_1/n_2=100$]{\includegraphics[width=0.45\textwidth]{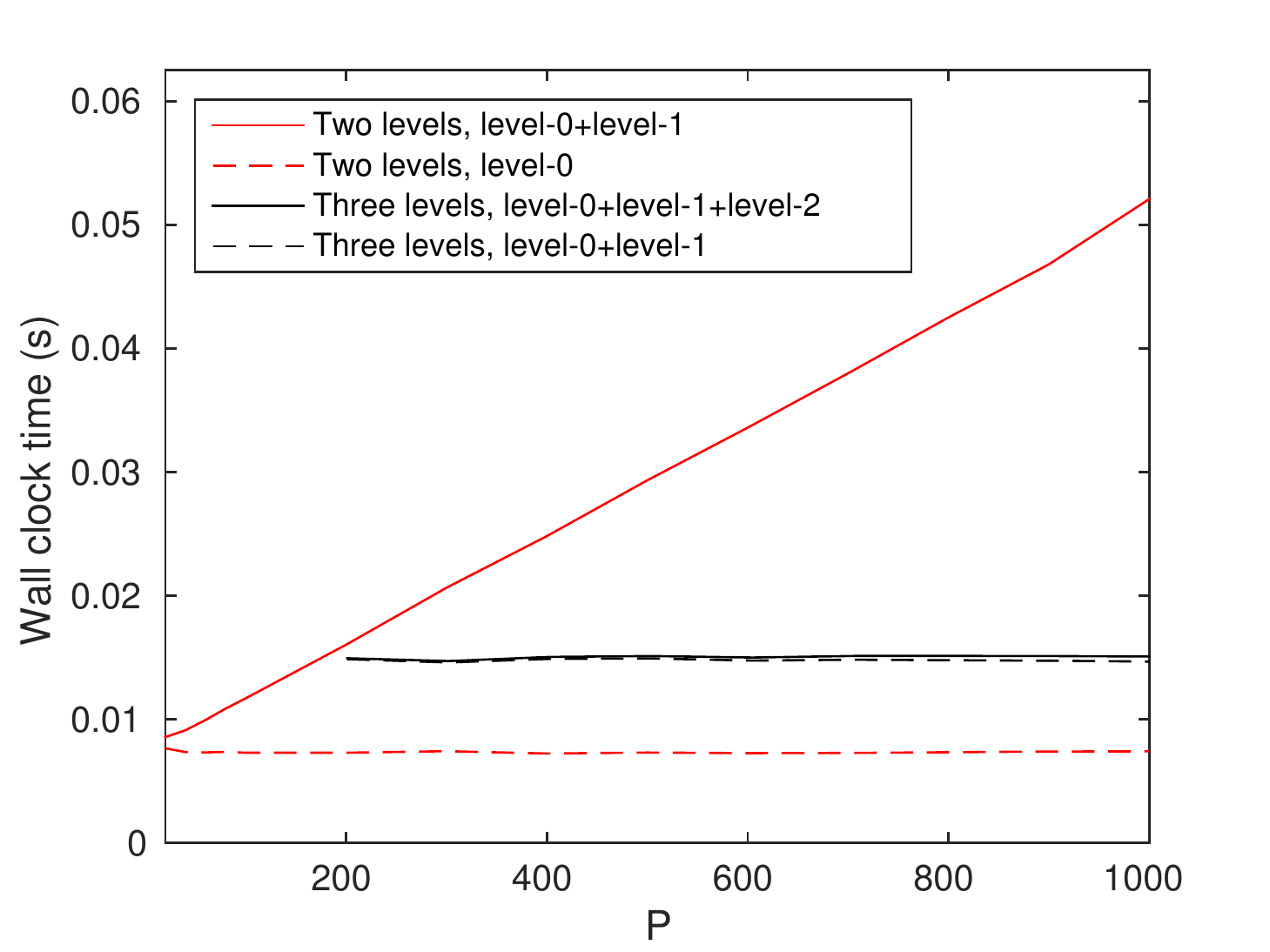}\label{fig-l2_l3_100}   }%
\caption{Comparison between the two-level ($\ell=1$) and the three-level ($\ell=2$) Schur complement solver with different local sizes. Times for local solvers in the level-0, local solvers in the level-1 and the Schur complement solve in level-2.} 
\label{fig-l2_l3}
\end{figure} 

In Fig. \ref{fig-l3fa}, an adaptive coarsening at level-2 is presented. Since $n_2 < \frac{n_1}{n_2}$ in the plot, we can reduce the coarsening from level-1 to level-2 to better balance the problems at these levels. It implies to enforce that $n_2 = \frac{n_1}{n_2}$. As expected, this choice is more efficient compared to the fixed size approach, when the size of the level-2 problem is below level-1 local problem sizes. 

 \begin{figure}[h!!]
\centering
\subfigure[$n_0/n_1=20$]{\includegraphics[width=0.42\textwidth]{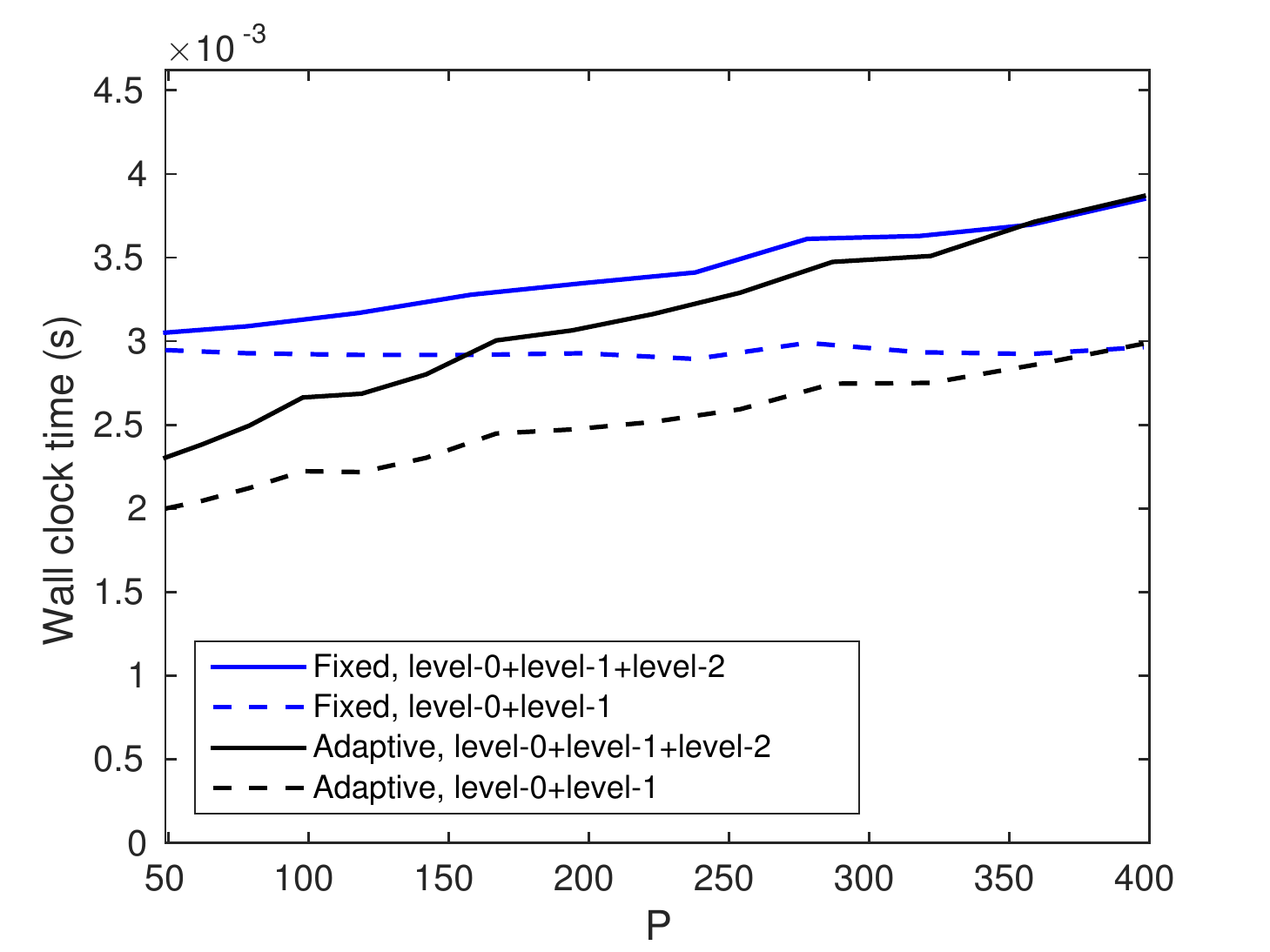} 
\label{fig-l3fa_20}}
\subfigure[$n_0/n_1=50$]{\includegraphics[width=0.42\textwidth]{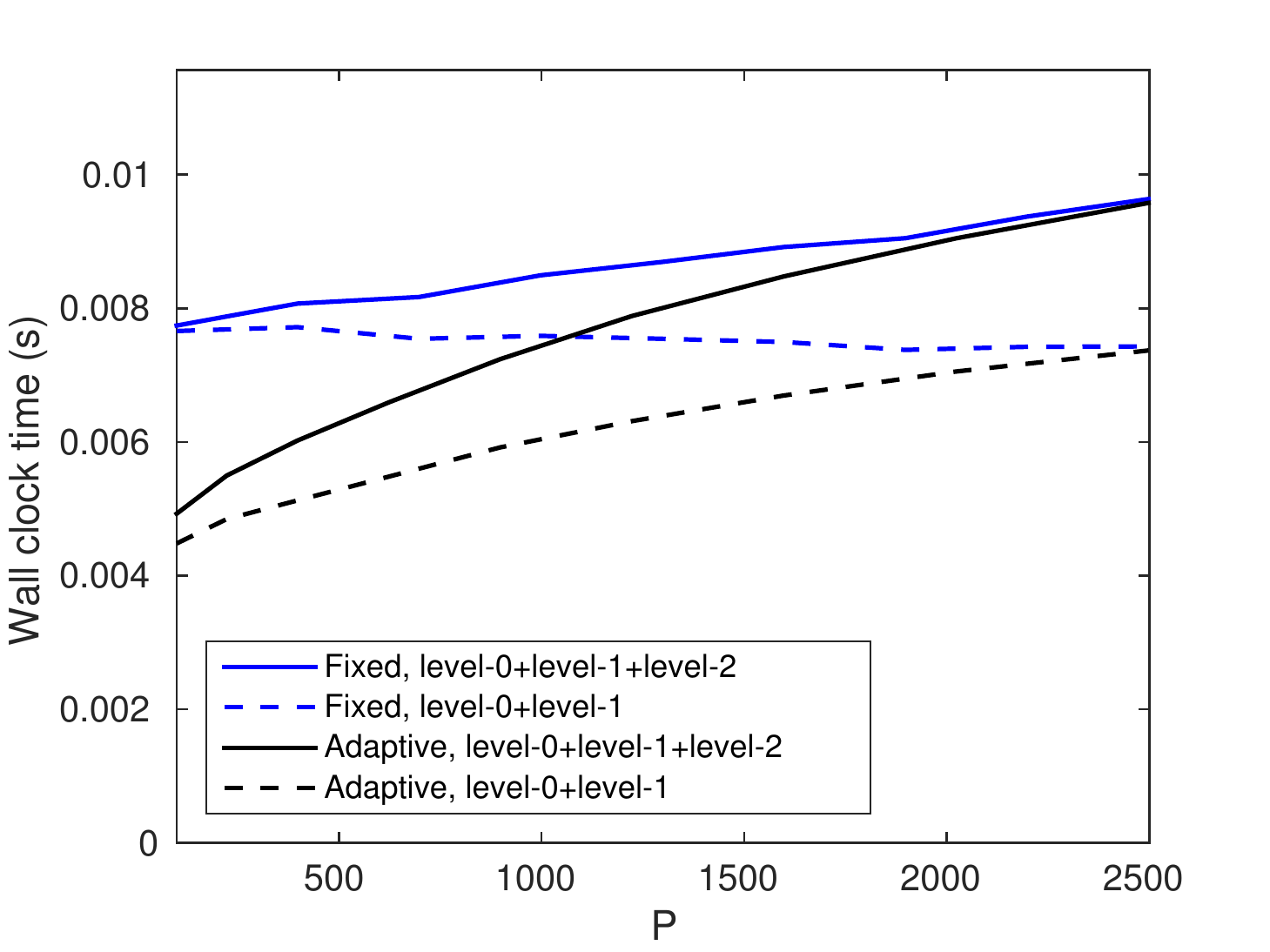} 
\label{fig-l3fa_50}}
\caption{ Comparison between the fixed coarsening and the adaptive coarsening between level-1 and level-2 for different level-0 local problem sizes. The adaptive coarsening uses $n_2=n_1/n_2=\sqrt{n_1}$. } 
\label{fig-l3fa}
\end{figure} 

\section{Conclusions}\label{sec:conclusions}

In this work, we have considered a parallel-in-time multilevel direct method and two iterative parallel-in-time nonlinear solvers for the numerical approximation of ODEs using parallel computations. The time-parallel direct method computes explicitly multi-level Schur complements. It allows for arbitrary high levels of concurrency and has very good theoretical speed-up ratios compared to traditional parareal methods. The method can be considered a parareal method with an automatic definition of the coarse solver. However, such definition makes the method to converge in one iteration, i.e., it is a direct method, and thus does not suffer from the poor parallel eficiency of parareal schemes. The proposed scheme can be applied to $\theta$-methods, DG methods, Runge-Kutta methods, and BDF methods. Next, we consider nonlinear ODEs, and propose two different strategies. First, we consider a Newton method over the global nonlinear ODE, using the multilevel Schur complement solver at every nonlinear iteration. Second, we state the global nonlinear problem in terms of the nonlinear Schur complement (at an arbitrary level), and perform nonlinear iterations over it. Such approach leads to nested nonlinear multilevel iterations.  

A detailed set of numerical experiments has been considered. Out of these results, the proposed methodology is observed to exhibit excellent scalability properties. The methods are weakly scalable in time, i.e., increasing X times the number of MPI tasks one can solve X times more time steps, in approximately the same amount of time, which is a key property to reduce time-to-solution in ODE simulations with heavy time stepping.

We refer to \cite{badia_space-time_2017} for the combination of these ideas with space-parallel highly scalable BDDC implementations \cite{badia_highly_2014,badia_multilevel_2016,badia_scalability_2015} to design space-time preconditioners for transient PDEs. The nonlinear algorithms in \cite{badia_space-time_2017} only consider the Newton-Schur complement solver. Future extensions of this work is the application of the level-$k$ nonlinear Schur complement scheme in space, leading to new nonlinear space-time preconditioners.

\bibliographystyle{siam_no_dashed}
\bibliography{art023}

\end{document}

%% file: fig_time_part.tex
  \begin{tikzpicture}
  \draw[arrows=->,line width=1pt]((-0.5,0)--(16.5,0);
  \foreach \i in{0,...,32}
  {
    \draw[fill=red] (0.5*\i,0) circle (0.15em);
  }
  \foreach \i in{0,1,2}
  {
    \draw[fill=black] (0.5*\i,0) 
    node[above] { $t_0^{\i}$};
  }
    \draw[fill=black] (1.5,0.1) 
    node[above] { $\cdots$};
    \draw[fill=black] (16,0)
    node[above] { $t_0^{n_0}$};

    \draw[arrows=->,line width=1pt]((-0.5,-1)--(16.5,-1);
    \foreach \i in{0,...,32}
  {
    \draw[fill=black] (0.5*\i,-1) circle (0.15em);
  }
  \foreach \i in{0,1,2}
  {
    \draw[fill=black] (2*\i+0.2,-1) 
    node[above] { $t_0^{m(\i)}$};
  }
  \draw[fill=black] (16+0.2,-1) 
    node[above] { $t_0^{m(n_1)}$};
  
  \foreach \i in {0,...,8}
  {
    \draw[fill=red] (2*\i,-1) circle (0.25em);
  }
  \foreach \i in {0,1,2}
  {
    \draw[fill=black] (2*\i,-1) 
    node[below] { $t_1^{\i}$};
  }
    \draw[fill=black] (6,-1.1) 
    node[below] { $\cdots$};
    \draw[fill=black] (16,-1) 
    node[below] { $t_1^{n_1}$};
    
    \draw[fill=black] (0,-1.8)
    node {$\vdots$};

    \foreach \i in{0,...,32}
  {
    \draw[fill=black] (0.5*\i,-3) circle (0.15em);
  }
  \foreach \i in {0,...,8}
  {
    \draw[fill=black] (2*\i,-3) circle (0.25em);
  }
  
  \draw[arrows=->,line width=1pt]((-0.5,-3)--(16.5,-3);
  \foreach \i in {0,8,16}
  {
    \draw[fill=red] (\i,-3) circle (0.35em);
  }
  \draw[fill=black] (0,-3.1) 
  node[below] { $t_\ell^0$};
  \draw[fill=black] (8,-3.2) 
  node[below] { $\cdots$};
  \draw[fill=black] (16,-3.1) 
  node[below] { $t_\ell^{n_\ell}$};
  \draw[fill=black] (0,-2.9) 
  node[above] { $t_{\ell-1}^{m(0)}$};
  \draw[fill=black] (16,-2.9) 
  node[above] { $t_{\ell-1}^{m(n_\ell)}$};
  

\end{tikzpicture}